\documentclass[leqno, myheadings, twoside]{amsart}

\usepackage{amsmath, amsthm, amssymb, amscd, amsxtra,graphicx}
\usepackage{latexsym, amsfonts}
\usepackage{easybmat}
\usepackage{etex}
\usepackage{url}
\usepackage{texdraw}
\usepackage{epsfig}
\usepackage{tikz}
\allowdisplaybreaks[4]
\usepackage{xcolor, soul}

\usepackage{geometry}
\usepackage{geometry}
\makeatletter \@addtoreset{equation}{section}

\makeatletter \renewcommand{\@biblabel}[1]{#1.}
\theoremstyle{remark}

\begin{document}
\setcounter{page}{1}
\title[Geometric invariants of spectrum of the Navier-Lam\'{e} operator]{Geometric invariants of spectrum of the Navier-Lam\'{e} operator}

\author{Genqian Liu}
\address{School of Mathematics and Statistics, Beijing Institute of Technology, Beijing 100081, China}
\email{liugqz@bit.edu.cn}
\subjclass[2010]{}
%\date={}
%\date{ }
%\thanks { }

\maketitle

\date{}
\protect
\maketitle ~~~\\[-15mm]

\begin{center}
{\footnotesize   School of Mathematics and Statistics, Beijing Institute of Technology, Beijing 100081, China\\
 Emails:  liugqz@bit.edu.cn \\
 }
\end{center}

%\\[5mm]

\vskip 0.59 true cm

\begin{abstract}    For a compact connected Riemannian $n$-manifold $(\Omega,g)$ with smooth boundary, we explicitly calculate the first two coefficients $a_0$ and $a_1$ of the asymptotic expansion of $\sum_{k=1}^\infty e^{-t \tau_k^\mp}= a_0t^{-n/2} \mp  a_1 t^{-(n-1)/2}+a_2^\mp t^{-(n-2)/2} +\cdots+ a_m^\mp t^{-(n-m)/2} +O(t^{-(n-m-1)/2})$ as $t\to 0^+$, where $\tau^-_k$ (respectively, $\tau^+_k$) is the $k$-th Navier-Lam\'{e} eigenvalue on $\Omega$ with Dirichlet (respectively, Neumann) boundary condition. These two coefficients provide precise information for the volume  of the elastic body $\Omega$ and the surface area of the boundary $\partial \Omega$ in terms of the spectrum of the  Navier-Lam\'{e} operator.
This gives an answer to an interesting and open problem mentioned by Avramidi in \cite{Avr10}.  More importantly, our method is valid to explicitly calculate all the coefficients $a_l^\mp$, $2\le l\le m$, in the above asymptotic expansion.
As an application, we show that an $n$-dimensional ball is uniquely determined by its Navier-Lam\'{e} spectrum among all bounded elastic bodies with smooth boundary.  \end{abstract}

\vskip 1.69 true cm

\section{ Introduction}

\vskip 0.45 true cm

For the Navier-Lam\'{e} elastic wave equations, one of the most important problems is to study the shape of the elastic body from its
 vibrational frequencies, because this kind of geometric property reveals the essential behavior of the elastic body.

Let $(\Omega,g)$ be a Riemannian $n$-manifold with smooth boundary
 $\partial \Omega$. Let $P_g$ be the Navier-Lam\'{e} operator:
 \begin{eqnarray} \label {1-1} P_g\mathbf{u}:=\mu \nabla^* \nabla \mathbf{u} -(\mu +\lambda) \,\mbox{grad}\; \mbox{div}\, \mathbf{u} -\mu \, \mbox{Ric} (\mathbf{u}), \;  \;\; \mathbf{u}=(u^1, \cdots, u^n),\end{eqnarray}
 where  $\mu$ and $\lambda$ are Lam\'{e} parameters satisfying $\mu>0$ and $\mu+\lambda\ge 0$,  $\nabla^* \nabla $ is the Bochner Laplacian (see (\ref{19.9.30-1}) in section 2, or (2.12) of \cite{Liu-1}), $\mbox{div}$ and $\mbox{grad}$ are the usual divergence and gradient operators, and  \begin{eqnarray} \label{18/12/22} \mbox{Ric} (\mathbf{u})= \big(\sum\limits_{k,l=1}^n R^{k\,\,1}_{\,lk} u^l,  \sum\limits_{k,l=1}^n R^{k\,\,2}_{\,lk} u^l, \cdots, \sum\limits_{k,l=1}^n R^{k\,\,n}_{\,lk} u^l\big)\end{eqnarray} denotes the action of Ricci tensor $\mbox{R}_l^{\;j}:=\sum_{k=1}^n R^{k\,\,j}_{\,lk}$ on $\mathbf{u}$.

   The natural boundary conditions for $P_g$ (see, Lemma 2.1.1 in \cite{Liu-1}) include prescribing $\mathbf{u}\big|_{\partial \Omega}$, Dirichlet type, and
  \begin{eqnarray}\label{21/6/9;1}  \frac{\boldsymbol{\partial} \mathbf{u}}{\partial \nu}:= 2\mu \,(\mbox{Def}\;\mathbf{u})^\# {\nu} +\lambda (\mbox{div}\; \mathbf{u}){\nu} \quad \mbox{on}\; \, \partial \Omega, \end{eqnarray}
   Neumann type, where  $(\mbox{Def}\, \mathbf{u})_{jk}= \frac{1}{2} \big(u_{j;k} +u_{k;j}\big)$ is the strain tensor, $\#$ is the sharp operator (for a tensor) by raising index  and $\nu$ is the unit inner normal to $\partial \Omega$. We denote by $P_g^-$ and $P_g^+$ the Navier-Lam\'{e} operators with the Dirichlet and Neumann boundary conditions, respectively.
    Since $P_g^-$ (respectively, $P_g^+$) is an unbounded, self-adjoint and positive (respectively,  nonnegative) operator in $[H^1_0(\Omega)]^n$ (respectively, $[H^1(\Omega)]^n$) with discrete spectrum $0< \tau_1^- < \tau_2^- \le \cdots \le \tau_k^- \le \cdots \to +\infty$ (respectively, $0\le \tau_1^+ < \tau_2^+ \le \cdots \le \tau_k^+ \le \cdots \to +\infty$),
 one has (see \cite{Avr10}, \cite{Pleij}, \cite{Agm}, \cite{Liu-1}, \cite{Lapt} or \cite{Hook})
  \begin{eqnarray} \label{1-4} P_g^\mp {\mathbf{u}}_k^\mp =\tau_k^\mp {\mathbf{u}}_k^\mp,\end{eqnarray} where ${\mathbf{u}}_k^-\in  [H^1_0(\Omega)]^n$ (respectively, ${\mathbf{u}}_k^+\in  [H^1(\Omega)]^n$) is the eigenvector corresponding to eigenvalue $\tau_k^{-}$ (respectively, $\tau_k^{+}$).
(\ref{1-4}) can be rewritten as
  \begin{eqnarray} \label{1-5} \left\{\! \begin{array}{ll}   \mu \nabla^*\nabla {\mathbf{u}}_k^- - (\mu+\lambda) \,\mbox{grad}\; \mbox{div}\;  {\mathbf{u}}_k^- -\mu\, \mbox{Ric} \,({\mathbf{u}}_k^-) = \tau^-_k {\mathbf{u}}_k^-  \;\; &\mbox{in}\;\;  \Omega,\\
        {\mathbf{u}}_k^-=0 \;\;& \mbox{on}\;\;  \partial \Omega
    \end{array}  \right.\end{eqnarray}
    and
    \begin{eqnarray} \label{1-6} \left\{ \!\begin{array}{ll}  \mu \nabla^*\nabla {\mathbf{u}}_k^+ - (\mu+\lambda)\, \mbox{grad}\; \mbox{div}\;  {\mathbf{u}}_k^+ -\mu\, \mbox{Ric} \,({\mathbf{u}}_k^+)   = \tau^+_k {\mathbf{u}}_k^+  \;\; &\mbox{in}\;\;  \Omega,\\
       \frac{\boldsymbol{\partial}  {\mathbf{u}}_k^+}{\partial \nu}=0 \;\;& \mbox{on}\;\;  \partial \Omega.
    \end{array}  \right.\end{eqnarray}

    Clearly, the eigenvalue problems (\ref{1-5}) and (\ref{1-6}) can be immediately obtained by considering the solutions of the form ${\mathbf{v}}(t, x)=T(t) {\mathbf{u}}(x)$ in the
    following  Navier-Lam\'{e} elastodynamic wave equations:
    \begin{eqnarray} \label{1-.7} \;\;\quad\;\;\; \left\{\! \begin{array}{ll}  \frac{\partial^2 {\mathbf{v}}^-}{\partial t^2}+ \mu \nabla^*\nabla {\mathbf{v}}^- - (\mu+\lambda) \,\mbox{grad}\; \mbox{div}\;  {\mathbf{v}}^- -\mu\, \mbox{Ric} \,({\mathbf{v}}^-)  =0
           \;\; &\mbox{in}\;\;  (0,+\infty)\times \Omega,\\
            {\mathbf{v}}^-=0 \;\;& \mbox{on}\;\; (0,+\infty)\times \partial \Omega, \\
     {\mathbf{v}}^-(0,x)= \mathbf{v}_0 \;\quad \frac{\partial {\mathbf{v}}^-}{\partial t}(0,x)=0 & \mbox{on}\;\; \{0\}\times \Omega \end{array}  \right.\end{eqnarray}
       and   \begin{eqnarray} \label{1-.8} \;\;\quad\;\;\; \left\{\! \begin{array}{ll}  \frac{\partial^2 {\mathbf{v}}^+}{\partial t^2}+ \mu \nabla^*\nabla {\mathbf{v}}^+ - (\mu+\lambda) \,\mbox{grad}\; \mbox{div}\;  {\mathbf{v}}^+ -\mu\, \mbox{Ric} \,({\mathbf{v}}^+)  =0
           \;\; &\mbox{in}\;\;  (0,+\infty)\times \Omega,\\
            \frac{\boldsymbol{\partial} {\mathbf{v}}^+}{\partial \nu}=0 \;\;& \mbox{on}\;\; (0,+\infty)\times \partial \Omega, \\
     {\mathbf{v}}^+(0,x)= \mathbf{v}_0 \;\quad \frac{\partial {\mathbf{v}}^+}{\partial t}(0,x)=0 & \mbox{on}\;\; \{0\}\times \Omega. \end{array}  \right.\end{eqnarray}
   In three spatial dimensions, the elastic wave equations (\ref{1-.7}) and (\ref{1-.8}) describe the propagations of waves in an isotropic homogeneous elastic medium. The elasticity of the material provides the restoring force of the wave. Most solid materials are elastic, so these two equations describe such phenomena as seismic waves in the Earth, ultrasonic waves used to detect flaws in materials and the deformations of thin elastic shells whose middle surface must stay inside a given surface in the three-dimensional Euclidean space.
  Equivalently, we can also rewrite the equation in (\ref{1-.7}) (or (\ref{1-.8})) into another form (cf. \cite{Liu-1} and \cite{Liu-2}), as it must account for both longitudinal and transverse motion (in three spatial dimensions):
 \begin{eqnarray}\label{2021.2.6-1} \frac{\partial^2 {\mathbf{v}}^\mp}{\partial t^2} +
          \mu \,\mbox{curl}\; \mbox{curl}\; {\mathbf{v}}^\mp  -(2\mu +\lambda) \,\mbox{grad}\; \mbox{div}\; {\mathbf{v}}^\mp -2\mu \,\mbox{Ric}\; ({\mathbf{v}}^{\mp})=0.\end{eqnarray}
 In particular, if $\mbox{div}\, {\mathbf{v}}^\mp$ are set to zero, (\ref{2021.2.6-1}) becomes (effectively) Maxwell's equations for the propagation of the electric field $\mathbf{v}^\mp$, which has only transverse waves (see \cite{Grif}, \cite{Coo}, \cite{Stra}). In addition, if $\lambda+\mu=0$, then (\ref{1-.7}) and (\ref{1-.8}) reduce to the classical wave equations, which has only longitudinal waves. For the derivation of the Navier-Lam\'{e} elastic wave equations, its mechanical meaning and the explanation of the Dirichlet and Neumann boundary conditions, we refer the reader to \cite{Liu-1} for the case of Riemannian manifold and to \cite{Bant}, \cite{Ev}, \cite{Gur}, \cite{Kaw}, \cite{LaLi}, \cite{KLV}, \cite{Sl} for the case of Euclidean space.

 The Navier-Lam\'{e} eigenvalues are physical quantities because they just are the square of vibrational frequencies of an elastic body in two or three dimensions. And these basic physical quantities can be measured experimentally. An interesting question, which is similar to the famous Kac question for the Dirichlet-Laplacian (see \cite{Kac}, \cite{Lo} or \cite{We1} ), is: ``can one hear the shape of an elastic body by hearing the vibrational frequencies (or pitches) of the elastic body?''
More precisely, does the spectrum of the Navier-Lam\'{e} operator determine the geometry of an elastic body (see \cite{Avr10})?

In the special case of $\mu+\lambda=0$ (i.e., the Navier-Lam\'{e} operator reduces to the Laplace operator), a celebrated result of the spectral
 (geometric) invariants had been obtained by McKean and Singer \cite{MS}. They proved the famous Kac conjecture and gave an explicit expression to the first three coefficients of asymptotic expansion for the  heat trace of the Laplacian on a bounded domain $\Omega$ of a Riemannian manifold:
 \begin{eqnarray} \label{18-03-20-1} \sum_{k=1}^\infty e^{-\beta_k^\mp t}\!\!\!&\!=\!\!\!&(4\pi t)^{-n/2}\Big(\mbox{Vol}(\Omega) \mp \frac{1}{4} \sqrt{4\pi t} \, \mbox{Vol}(\partial \Omega)\\
 && +\frac{t}{3}\, \int_{\Omega} R -\frac{t}{6} \, \int_{\partial \Omega}J +O(t^{3/2})\Big) \;\;\,\,\mbox{as}\,\, t\to 0^+,\nonumber \end{eqnarray}
where $\beta_k^{-}$ (respectively, $\beta_k^{+}$) is the $k$-th Dirichlet-Laplacian (respectively, Neumann-Laplacian) eigenvalue on $\Omega$; $R$ and $J$ are the scalar curvature and the mean curvature of $\Omega$ and $\partial \Omega$, respectively; $O(t^{3/2})$ cannot be improved.

 The symbolic approach by Seeley \cite{See,See2} and Greiner \cite{Gre} is a very powerful general analytical
procedure for analyzing the structure of the asymptotic expansion based on
the theory of pseudodifferential operators and the calculus of symbols of operators. This approach may be considered
for calculation of the heat invariants explicitly in terms of the jets of the symbol
of the operator; it provides an iterative procedure for such a calculation. However,
as far as we know, because of the technical complexity and, most importantly,
lack of the manifest covariance, such analytical tools have never been used for the
actual calculation of the explicit form of the heat invariants in an invariant geometric
form (see \cite{Avr10}).
 The systematic explicit calculation of heat kernel coefficients for Laplace type operators is now well understood due to the work of Gilkey
\cite{Gil} and many others (see \cite{Gil2, Gilk, Gr, Gr2, Vas, Kir, BGV, Avr, LiuA, AGMT} and references therein) because the Riemannian structure on a manifold is determined by a Laplace type operator. For the classical boundary conditions, like Dirichlet, Neumann, Robin,
and mixed combination thereof on vector bundles, the coefficients of the trace of  heat kernel have been
explicitly computed up to the first five terms (see, for example, \cite{Kir4, BG, BGKV, Avra}). For other type operators (which originated from physics problems), the corresponding explicit form of the heat invariants have also been discussed.
    Liu in \cite{Liu1} explicitly calculated
    the first two coefficients of asymptotic expansion of the heat trace for the Stokes operator (i.e., incompressible slow flow operator),
   and in \cite{Liu2, Liu3} gave the first four coefficients of asymptotic expansion of the heat trace for the Dirichlet-to-Neumann map (We also refer the reader to \cite{PolS} for the asymptotic expansion of the first three coefficients) as well as polyharmonic Steklov operator.

Contrary to the Laplace type operators, there are no systematic effective methods
for an explicit calculation of the spectral invariants for second-order operators
which are not of Laplace type. Such operators appear in so-called matrix geometry
\cite{Avr7, Avr8, Avr9}, when instead of a single Riemannian metric there is a matrix-valued
symmetric $2$-tensor. Let us point out that the Navier-Lam\'{e} operator is just a non-Laplace type operator. Five decades ago Greiner \cite{Gre}, p.$\,$164, indicated that ``the problem of interpreting
these coefficients geometrically remains open''. There has not been much
progress in this direction. Thus, {\it the geometric aspect of
the spectral asymptotics of the Navier-Lam\'{e} operator remains an open problem} (see \cite{Avr10}).
In the celebrated paper \cite{Avr10}, Avramidi consider general non-Laplace type operators on manifolds with boundary by introducing a ``noncommutative''
Dirac operator as a first-order elliptic partial differential operator
such that its square is a second-order self-adjoint elliptic operator with positive
definite leading symbol (not necessarily of Laplace type) and study the spectral
asymptotics of these operators with Dirichlet boundary conditions. However, Avramidi's method is too complicated to be applied to the Navier-Lam\'{e} operator because of extremely technical difficulty to calculate the noncommutative Dirac operator and required symbol integral for such an operator.

In this paper, by combining the technique of calculus of symbols for the integral kernel and  ``method of images'', we obtain the following result:

\vskip 0.25 true cm

  \noindent{\bf Theorem 1.1.} \ {\it Let $(\Omega,g)$ be a compact Riemannian manifold of dimension $n$ with smooth boundary $\partial \Omega$, and let $0< \tau_1^-< \tau_2^- \le \cdots \le \tau_k^- \le \cdots$
 (respectively, $0\le \tau_1^+ < \tau_2^+ \le \tau_3^+ \le \cdots \le \tau_k^+ \le \cdots $) be the eigenvalues
 of the Navier-Lam\'{e} operator $P_g^-$ (respectively, $P_g^+$) with respect to the Dirichlet (respectively, Neumann) boundary condition.  Then
 \begin{eqnarray} \label{1-7} &&  \sum_{k=1}^\infty e^{-\tau_k^\mp t}  = \mbox{Tr}\,(e^{-tP_g^\mp})=\bigg[ \frac{(n-1)}{(4\pi \mu t)^{n/2}}
 + \frac{1}{(4\pi (2\mu+\lambda) t)^{n/2}}\bigg] {\mbox{Vol}}\,(\Omega) \\
&&  \quad \, \mp \frac{1}{4} \bigg[  \frac{(n-1)}{(4\pi \mu t)^{(n-1)/2}}
 +  \frac{1}{(4\pi (2\mu+\lambda) t)^{(n-1)/2}}\bigg]{\mbox{Vol}}\,(\partial\Omega)+O(t^{{1-n}/2})\quad\;\; \mbox{as}\;\; t\to 0^+.\nonumber\end{eqnarray}
 Here ${\mbox{Vol}}\,(\Omega)$ denotes the $n$-dimensional volume of $\Omega$,  ${\mbox{Vol}}\,(\partial\Omega)$ denotes the $(n-1)$-dimensional volume of $\partial \Omega$.}

 \vskip 0.25 true cm

  Our result shows that not only the volume ${\mbox{Vol}}(\Omega)$  but also the surface area $\mbox{Vol}(\partial \Omega)$ can be obtained if we know all the Navier-Lam\'{e} eigenvalues with respect to the Dirichlet (respectively, Neumann) boundary condition. This gives an answer to an interesting and open problem mentioned by Avramidi in \cite{Avr10}. Roughly speaking, one can ``hear'' the volume of the domain and the surface area of its boundary by ``hearing'' all the pitches of the vibration of an elastic body.

\vskip 0.28true cm

The key ideas of this paper are as follows. We denote by $(e^{-tP_g^\mp})_{t\ge0}$  the  parabolic semigroups generated by $-P_g^\mp$.
 More precisely, ${\mathbf{w}}^\mp (t,x) = e^{-tP_g^\mp} {\mathbf{w}}_0(x)$ solve the following initial-boundary problems: \begin{eqnarray} \label{1-2}\quad\; \quad\;\; \left\{ \!\! \begin{array}{ll}  \frac{\partial {\mathbf{w}}^-}{\partial t} + \mu \nabla^*\nabla {\mathbf{w}}^- - (\mu+\lambda) \,\mbox{grad}\; \mbox{div}\;  {\mathbf{w}}^- -\mu\, \mbox{Ric} \,({\mathbf{w}}^-) =0
           \; &\mbox{in}\;\;  (0,+\infty)\times \Omega,\\
            {\mathbf{w}}^-=0 \;\;& \mbox{on}\;\; (0,+\infty)\times \partial \Omega, \\
     {\mathbf{w}}^-(0,x)= \mathbf{w}_0 \;\; & \mbox{on}\;\; \{0\}\times \Omega \end{array}  \right.\end{eqnarray}
       and   \begin{eqnarray} \label{1-3} \quad \; \quad\;\;\left\{ \!\!\begin{array}{ll}  \frac{\partial {\mathbf{w}}^+}{\partial t} + \mu \nabla^*\nabla {\mathbf{w}}^+ - (\mu+\lambda) \,\mbox{grad}\; \mbox{div}\;  {\mathbf{w}}^+ -\mu\, \mbox{Ric} \,({\mathbf{w}}^+) =0
           \; &\mbox{in}\;\;  (0,+\infty)\times \Omega,\\
            \frac{\boldsymbol{\partial} {\mathbf{w}}^+}{\partial \nu}=0 \;\;& \mbox{on}\;\; (0,+\infty)\times \partial \Omega, \\
     {\mathbf{w}}^+(0,x)= \mathbf{w}_0 \;\; & \mbox{on}\;\; \{0\}\times \Omega. \end{array}  \right.\end{eqnarray}
   If $\{{\mathbf{u}}_k^\mp\}_{k=1}^\infty$ are orthonormal eigenvectors of the Navier-Lam\'{e} problem corresponding to eigenvalues $\{\tau_k^\mp\}_{k=1}^\infty$, then the integral kernels  ${\mathbf{K}}^\mp(t, x, y)=e^{-t P_g^\mp} \delta(x-y)$ of the semigroups  are given by \begin{eqnarray} \label{1-0a-1} {\mathbf{K}}^\mp(t,x,y) =\sum_{k=1}^\infty e^{-t \tau_k^\mp} {\mathbf{u}}_k^\mp(x)\otimes {\mathbf{u}}_k^\mp(y).\end{eqnarray}
Thus the integrals of the traces of ${\mathbf{K}}^\mp(t,x,y)$ are actually spectral invariants:
\begin{eqnarray} \label{1-0a-2}\int_\Omega  \big(\mbox{Tr}\,({\mathbf{K}}^\mp(t,x,x))\big) dV=\sum_{k=1}^\infty e^{-t \tau_k^\mp}.\end{eqnarray}
To further analyze the geometric content of the spectrum, we calculate the same
 integral of the trace by another approach:
  let $\mathcal{M}=\Omega\cup (\partial \Omega) \cup \Omega^*$ be the (closed) double of $\Omega$, and $\mathcal{P}$ the double to $\mathcal{M}$ of the operator $P_g$ on $\Omega$ (see section 4). Then $-\mathcal{P}$ generates a strongly continuous semigroup $(e^{-t\mathcal{P}})_{t\ge 0}$ on $L^2(\mathcal{M})$ with integral kernel
  $\mathbf{K} (t, x,y)$.  Clearly, ${\mathbf{K}}^{\mp} (t,x,y)= {\mathbf{K}}(t,x,y)  \mp {\mathbf{K}} (t,x,\overset{*}{y})$ for $x,y\in \bar \Omega$,
  where $\overset{*} {y}$ is the double of $y\in \Omega$. This technique stems from  McKean and Singer (see \cite{MS}), and is called ``method of images''.
  Since $e^{-t\mathcal{P}}f(x) =\frac{1}{2\pi i} \int_{\mathcal{C}} e^{-t\tau } (\tau I- \mathcal{P})^{-1} f(x) \, d\tau$, we have
  $$e^{-t\mathcal{P}} f(x) = \frac{1}{(2\pi)^n} \int_{{\mathbb{R}}^n} e^{ix\cdot \xi} \Big( \frac{1}{2\pi i} \int_{\mathcal{C}} e^{-t\tau}\, \iota\big((\tau I-\mathcal{P})^{-1}\big) \, \hat{f} (\xi) \, d\tau\Big) d\xi,$$
  so that \begin{eqnarray*}&&\mathbf{K}(t,x,y)= e^{-t\mathcal{P}} \delta(x-y) = \frac{1}{(2\pi)^n} \int_{{\mathbb{R}}^n} e^{i(x-y)\cdot \xi} \Big( \frac{1}{2\pi i} \int_{\mathcal{C}}  e^{-t\tau}\, \iota\big((\tau I-\mathcal{P})^{-1}\big) \, d\tau\Big) d\xi\\
  && \qquad \qquad  \; = \frac{1}{(2\pi)^n} \int_{{\mathbb{R}}^n} e^{i(x-y)\cdot \xi} \Big( \frac{1}{2\pi i} \int_{\mathcal{C}}  e^{-t\tau}\, \big(\sum_{l\ge 0} {\mathbf{q}}_{-2-l} (x, \xi,\tau)\big) \, d\tau\Big) d\xi,\end{eqnarray*}
  where $\mathcal{C}$ is a suitable curve in the complex plane in the positive direction around the spectrum of $\mathcal{P}$, and $\iota((\tau I-\mathcal{P})^{-1}):=\sum_{l\ge 0} {\mathbf{q}}_{-2-l} (x,\xi,\tau)$ is the full symbol of resolvent operator $(\tau I-\mathcal{P})^{-1}$.
  This implies that for any $x\in \bar \Omega$,
  \begin{eqnarray*}&&\mbox{Tr}\,( {\mathbf{K}}(t,x,x))= \frac{1}{(2\pi)^n} \int_{{\mathbb{R}}^n}   \Big( \frac{1}{2\pi i} \int_{\mathcal{C}}  e^{-t\tau}\,\sum_{l\ge 0}  \mbox{Tr}\,({\mathbf{q}}_{-2-l} (x,\xi, \tau) )\, d\tau\Big) d\xi,\\
  &&  \mbox{Tr}\,({\mathbf{K}}(t,x,\overset{*}{x}))= \frac{1}{(2\pi)^n} \int_{{\mathbb{R}}^n}  e^{i(x-\overset{*}{x})\cdot \xi}   \Big(  \frac{1}{2\pi i} \int_{\mathcal{C}}  e^{-t\tau}\,\sum_{l\ge 0} \mbox{Tr}\,({\mathbf{q}}_{-2-l} (x,\xi, \tau)) \, d\tau\Big) d\xi
  .\end{eqnarray*}
  It is easy to show that for any $x\in \Omega$,
   \begin{eqnarray*}&& \frac{1}{(2\pi)^n} \int_{{\mathbb{R}}^n}   \Big( \frac{1}{2\pi i} \int_{\mathcal{C}}  e^{-t\tau}\,\sum_{l\ge 1}  \mbox{Tr}\,({\mathbf{q}}_{-2-l} (x,\xi, \tau) )\, d\tau\Big) d\xi=O(t^{1-\frac{n}{2}}) \quad \mbox{as}\;\,t \to 0^+, \\
  &&  \frac{1}{(2\pi)^n} \int_{{\mathbb{R}}^n}  e^{i(x-\overset{*}{x})\cdot \xi}   \Big(  \frac{1}{2\pi i} \int_{\mathcal{C}}  e^{-t\tau}\,\sum_{l\ge 1} \mbox{Tr}\,({\mathbf{q}}_{-2-l} (x,\xi, \tau)) \, d\tau\Big) d\xi=O(t^{1-\frac{n}{2}})\quad \mbox{as}\;\, t\to 0^+.\end{eqnarray*}
    In order to finally establish asymptotic estimate, we denote by $U_\epsilon (\partial \Omega)$ the $\epsilon$-neighborhood of $\partial \Omega$ in $\mathcal{M}$. We can also show that
   \begin{eqnarray*} &&\ \frac{1}{(2\pi)^n} \int_{{\mathbb{R}}^n}   \Big( \frac{1}{2\pi i} \int_{\mathcal{C}}   {\mathbf{q}}_{-2} (x, \xi, \tau) \, d\tau\Big) d\xi= \left( \frac{n-1}{(4\pi \mu t)^{n/2}} + \frac{1}{ (4\pi (2\mu+\lambda)t)^{n/2}}\right) \quad \mbox{for}\;\, x\in \Omega,\\
      &&\; \,\frac{1}{(2\pi)^n} \int_{{\mathbb{R}}^n}  \! e^{i(x-\overset{*}{x})\cdot \xi}  \Big( \frac{1}{2\pi i} \int_{\mathcal{C}}    {\mathbf{q}}_{-2} (x,\xi, \tau) \, d\tau\Big) d\xi=O(t^{1-\frac{n}{2}})  \quad \mbox{as}\;\, t \to 0^+ \;\; \mbox{for}\;\, x\in \Omega\setminus U_\epsilon (\partial \Omega)   \end{eqnarray*} and
 \begin{eqnarray*} \!\!\!&&\int_{\Omega\cap U_\epsilon(\partial \Omega)}\!\left\{\!\frac{1}{(2\pi)^n}\! \int_{{\mathbb{R}}^n}  \! e^{i(x-\overset{*}{x})\cdot \xi}  \Big( \frac{1}{2\pi i} \!\int_{\mathcal{C}} \!    {\mathbf{q}}_{-2} (x,\xi, \tau) \, d\tau\Big) d\xi\!\right\}\!dV \\
 && \qquad \qquad \!=\! \frac{1}{4} \Big(\frac{n-1}{(4\pi \mu t)^{(n-1)/2}} \!+\! \frac{1}{(4\pi (2\mu\!+\!\lambda)t)^{(n-1)/2}}\Big) \mbox{Vol}\,(\partial \Omega)\!+\!O(t^{1-\frac{n}{2}}) \quad \mbox{as}\;\, t \to 0^+.
 \end{eqnarray*}
 Hence \begin{eqnarray} \int_{\Omega}\mbox{Tr}\big( {\mathbf{K}}^\mp(t,x,x) \big)dV = a_0t^{-n/2} \mp a_1 t^{-(n-1)/2} + O(t^{1-n/2})\quad \quad \mbox{as}\;\; t\to 0^+,\end{eqnarray}
where $a_0= \left(\frac{n-1}{(4\pi \mu )^{n/2}} + \frac{1}{ (4\pi (2\mu+\lambda))^{n/2}}\right){\mbox{Vol}}(\Omega)$ and $a_1=\frac{1}{4} \left(  \frac{(n-1)}{(4\pi \mu )^{(n-1)/2}}
 +  \frac{1}{(4\pi (2\mu+\lambda) )^{(n-1)/2}}\right){\mbox{Vol}}(\partial\Omega)$.

\vskip 0.18 true cm

As an application of theorem 1.1, we can prove the following spectral rigidity result:
\vskip 0.25 true cm

  \noindent{\bf  Theorem 1.2.} \ {\it Let $\Omega \subset {\mathbb{R}}^n$ be a bounded domain with smooth boundary $\partial \Omega$.
   Suppose that the Navier-Lam\'{e} spectrum with respect to the Dirichlet (respectively, Neumann) boundary condition, is equal to that of $B_r$, a ball of radius $r$. Then $\Omega=B_r$. }

\vskip 0.23 true cm

  Theorem 1.2 also shows that a ball is uniquely determined by its Navier-Lam\'{e} spectrum among all Euclidean bounded domains (elastic bodies) with smooth boundary.

\vskip 0.25 true cm

\vskip 1.49 true cm

\section{Some notations and lemmas}

\vskip 0.45 true cm

 Let $\Omega$ be an $n$-dimensional Riemannian  manifold (possibly with boundary), and let $\Omega$ be equipped with a smooth metric tensor $g =
\sum\limits_{j,k=1}^n g_{jk}\, dx_j \otimes dx_k$. Denote by $[g^{jk}]_{n\times n}$ the inverse of the matrix $[g_{jk}]_{n\times n}$ and set $|g|:= \mbox{det}\, [g_{jk}]_{n\times n}$. In
particular, $d\mbox{V}$, the volume element in $\Omega$ is locally given by $d\mbox{V} = \sqrt{|g|}\, dx_1\cdots dx_n$. By $T\Omega$ and $T^*\Omega$ we denote, respectively, the tangent and cotangent bundle on $\Omega$. Throughout, we shall also denote by $T\Omega$ global ($C^\infty$) sections in $T\Omega$ (i.e., $T\Omega \equiv C^\infty (\Omega, T\Omega)$);
similarly, $T^*\Omega \equiv C^\infty(\Omega, T^*\Omega)$.
A vector field on $\Omega$ is a section of the map $\pi: T\Omega\to\Omega$.  More concretely, a vector field is a smooth map $X: \Omega\to T\Omega$, usually written $p\mapsto X_p$, with the property that
\begin{eqnarray*} \pi \circ X={\mbox{Id}}_\Omega,\end{eqnarray*}
or equivalently, $X_p \in T_p\Omega$ for each $p\in \Omega$. If $(U; x_1,\cdots, x_n)$ is any smooth coordinate chart for $\Omega$, we can write the value of $X$ at any point $p\in U$ in terms of the coordinate basis vectors $\{\frac{\partial}{\partial x_j}\big|_{p}\}$ of $T_p\Omega$:
\begin{eqnarray*}  X_p= \sum_{j=1}^n X^j (p) \, \frac{ \partial}{\partial x_j}\bigg|_{p}.\end{eqnarray*}
This defines $n$ functions $X^j:U\to {\mathbb{ R}}$, called the component functions of $X$ in the given chart.
Recall first that
\begin{eqnarray} \label{9-2.1} \mbox{div} \,{X} :=\sum\limits_{j=1}^n \frac{1}{\sqrt{|g|}}\, \frac{\partial(\sqrt{|g|} \,{{X}}^j)}{\partial x_j}\quad \, \mbox{if}\;\; {{X}}=\sum\limits_{j=1}^n {{X}}^j \frac{\partial}{\partial x_j}\in T\Omega,\end{eqnarray}
and \begin{eqnarray} \label{9-2.2} \mbox{grad}\, v = \sum\limits_{j,k=1}^n \bigg(g^{jk} \frac{\partial v}{\partial x_k}\bigg)\frac{\partial}{\partial x_j}\quad \, \mbox{if}\;\; v\in C^\infty(\Omega), \end{eqnarray}
are, respectively, the usual divergence and gradient operators. Accordingly, the Laplace-Beltrami
operator $\Delta_g$ is just given by
 \begin{eqnarray} \label{18/12/22-2} \Delta_g:= \mbox{div}\; \mbox{grad} = \frac{1}{\sqrt{|g|}} \sum_{j,k=1}^n \frac{\partial}{\partial x_j} \bigg(\sqrt{|g|}\,g^{jk} \frac{\partial}{\partial x_k}\bigg). \end{eqnarray}
Next, let $\nabla$ be the associated Levi-Civita connection. For each ${X} \in T\Omega$, $\nabla X$ is the tensor of type $(0,2)$ defined by
\begin{eqnarray} \label{9-2.5} (\nabla {X})({Y},{Z}):= \langle \nabla_{{Z}} {X}, {Y}\rangle, \quad \; \forall\, {Y}, {Z}\in T\Omega. \end{eqnarray}
It is well-known that in a local coordinate system with the naturally associated frame field on the tangent bundle,
\begin{eqnarray*} \nabla_{\frac{\partial}{\partial x_k}} X = \sum\limits_{j=1}^n \big(\frac{\partial X^j}{\partial x_k} +\sum\limits_{l=1}^n \Gamma_{lk}^j X^l  \big)\frac{\partial }{\partial x_j}\quad \; \mbox{for}\;\; X=\sum\limits_{j=1}^n X^j \frac{\partial}{\partial x_j}, \end{eqnarray*}
 where $\Gamma_{lk}^j= \frac{1}{2} \sum_{m=1}^n g^{jm} \big( \frac{\partial g_{km}}{\partial x_l} +\frac{\partial g_{lm}}{\partial x_k} -\frac{\partial g_{lk}}{\partial x_m}\big)$ are the Christoffel symbols associated with the metric $g$ (see, for example, p.$\,$549 of \cite{Ta2}). If we denote \begin{eqnarray*} {X^j}_{;k}= \frac{\partial X^j}{\partial x_k} +\sum\limits_{l=1}^n \Gamma_{lk}^j X^l,\end{eqnarray*}
 then     \begin{eqnarray*} \label{18/10/29}  \nabla_Y X = \sum\limits_{j,k=1}^{n} Y^k {X^j}_{;k} \,\frac{\partial}{\partial x_j} \;\; \mbox{for}\,\; Y=\sum\limits_{k=1}^n Y^{k} \frac{\partial}{\partial x_k}.\end{eqnarray*}
  The symmetric part of $\nabla {X}$ is $\mbox{Def}\, {X}$, the deformation of ${X}$, i.e.,
\begin{eqnarray} (\mbox{Def}\; {X})({Y},{Z}) =\frac{1}{2} \{ \langle \nabla_{{Y}} {X}, {Z}\rangle +\langle \nabla_{{Z}} {X}, {Y}\rangle \}, \quad \, \forall\, {Y}, {Z}\in T\Omega\end{eqnarray}
(whereas the antisymmetric part of $\nabla\, {X}$ is simply $d{X}$, i.e.,
\begin{eqnarray*} d{X}({Y}, {Z}) = \frac{1}{2} \{ \langle \nabla_{{Y}} {X}, {Z}\rangle -\langle \nabla_{{Z}} {X}, {Y}\rangle \}, \quad \, \forall \,{Y}, {Z}\in T\Omega.)\end{eqnarray*}
Thus, $\mbox{Def}\; {X}$ is a symmetric tensor field of type $(0,2)$. In coordinate notation,
 \begin{eqnarray} \label{9-2.7} (\mbox{Def}\; {X})_{jk} = \frac{1}{2} ({{X}}_{j;k} +{{X}}_{k;j}), \quad \, \forall\, j,k,\end{eqnarray}
 where, ${{X}}_{k;j}:= \frac{\partial {{X}}_k}{\partial {\mathbf{x}}_j} -\sum_{l=1}^n \Gamma^{l}_{kj} {{X}}_l$ for a vector field ${X} =\sum_{j=1}^n {{X}}^j \frac{\partial}{\partial x_j}$, and $X_k=\sum_{l=1}^n g_{kl} X^l$.
 The adjoint  ${\mbox{Def}}^*$ of $\mbox{Def}$ is defined in local coordinates by $({\mbox{Def}}^* {w})^j=-\sum\limits_{k=1}^n {{w}}^{jk}_{\;\;\;\,;k}$ for each symmetric tensor field $w:=w_{jk}$ of type $(0,2)$.
The Riemann curvature tensor $\mathcal{R}$ of $\Omega$ is given by
\begin{eqnarray} \label{9-2.10} \mathcal{R}({X}, {Y}){Z} = [\nabla_{{X}}, \nabla_{{Y}}]{Z} -\nabla_{[{X},{Y}]} {Z}, \quad \, \forall\, {X}, {Y}, {Z} \in T\Omega, \end{eqnarray}
where $[{X}, {Y}] := {X}{Y} - {Y} {X}$ is the usual commutator bracket. It is convenient to change this
into a $(0, 4)$-tensor by setting
\begin{eqnarray*} \mathcal{R}({X}, {Y}, {Z}, {W}) := \langle \mathcal{R}({X}, {Y}){Z}, {W}\rangle, \quad \; \forall\, {X}, {Y}, {Z}, {W} \in T\Omega.
\end{eqnarray*}
In other words, in a local coordinate system such as that discussed above,
 \begin{eqnarray*} &&R_{jklm}= \bigg\langle \mathcal{R}\left(\frac{\partial}{\partial x_l}, \frac{\partial}{\partial x_m}\right) \frac{\partial}{\partial x_k}, \frac{\partial}{\partial x_j}\bigg\rangle.\end{eqnarray*}
The Ricci curvature $\mbox{Ric}$ on $\Omega$ is a $(0, 2)$-tensor defined as a contraction of $\mathcal{R}$:
\begin{eqnarray*} \mbox{Ric} ({X}, {Y} ):= \sum_{j=1}^n \bigg\langle \mathcal{R}\bigg(\frac{\partial }{\partial x_j}, {Y}\bigg) {X}, \frac{\partial }{\partial x_j}\bigg\rangle = \sum_{j=1}^n \bigg\langle \mathcal{R}\bigg({Y},\frac{\partial }{\partial x_j}\bigg)\frac{\partial }{\partial x_j}, {X}\bigg\rangle, \quad \forall\, {X}, {Y} \in T\Omega.\end{eqnarray*}
  That is,  \begin{eqnarray}\label{19.10.3-2}   R_{jk} = \sum\limits_{l=1}^n R^l_{jlk}=\sum\limits_{l,m=1}^n g^{lm}R_{ljmk}.\end{eqnarray}
 Note that \begin{eqnarray} \label{19.6.25-1} R^l_{jlk}=\frac{\partial \Gamma^l_{jk}}{\partial x_l}- \frac{\partial \Gamma^l_{jl}}{\partial x_k} +\sum_{s=1}^n\big(\Gamma^l_{sl} \Gamma^s_{jk}- \Gamma^l_{sk}\Gamma^s_{jl}\big).\end{eqnarray}

\vskip 0.17 true cm

In \cite{Liu-1}, by considering the equilibria states of elastic energy functional $\mathcal{E} (\mathbf{u})\!=\!-\frac{1}{2}\! \int_\Omega \!\big(\! \lambda (\mbox{div}\, \mathbf{u} )^2 \!+\!2 \mu \langle \mbox{Def}\; \mathbf{u}, \mbox{Def}\; \mathbf{u} \rangle \!\big) dV$, the author of this paper proved the following result, which generalizes the classical Navier-Lam\'{e} operator from the Euclidean space to a Riemannian manifold:

 \vskip 0.27 true cm

 \noindent{\bf Lemma 2.1.} \ {\it On a Riemannian manifold $(\Omega,g)$, modeling a homogeneous, linear, isotropic, elastic medium, the Navier-Lam\'{e} operator   ${\mathcal{L}}_g$ is given by}
  \begin{eqnarray} \label{18/11/4-4} P_g \mathbf{u} \!\!& \!\!\!=\!\!&\!\!\! \mu \, \nabla^* \nabla \mathbf{u} -(\mu+\lambda)\, \mbox{grad}\,\mbox{div}\, \mathbf{u}- \mu\,\mbox{Ric} (\mathbf{u}) \quad \mbox{for}\;\, \mathbf{u}=\sum\limits_{k=1}^n u^k\frac{\partial }{\partial x_k}\in T\Omega, \end{eqnarray}
 {\it where} $\nabla^*\nabla \mathbf{u}$ {\it is the Bochner Laplacian of} $\mathbf{u}$ {\it defined by}  \begin{align} \label{19.9.30-1} &\\
   & \nabla^*\nabla \mathbf{u}\!=\!- \sum\limits_{j=1}^n\!\bigg\{\! \Delta_g u^j  \!+\!2 \!\sum\limits_{\!m,k,l\!=\!1}^n\! \!g^{ml} \Gamma_{\!km}^j \frac{\partial u^k}{\partial x_l} \!+\!\!
   \sum_{\!m,k,l\!=\!1}^n \!\Big(\!g^{ml} \frac{\partial \Gamma^j_{\!kl}}{\partial x_m} \!+\!\sum\limits_{\!h\!=\!1}^n \!g^{ml} \Gamma_{\!hl}^j \Gamma_{\!km}^h \!-\!\sum\limits_{h=1}^n\! g^{ml} \Gamma_{\!kh}^j \Gamma_{ml}^h \!\Big)u^k\!\bigg\}\!\frac{\partial}{\partial x_j},\! \!\nonumber \end{align}
   and  \begin{eqnarray} \label{18/11/25} \mbox{Ric} (\mathbf{u})= \sum\limits_{j=1}^n  \big(\sum\limits_{k=1}^nR_k^j u^k\big)\frac{\partial }{\partial x_j}. \end{eqnarray}
  {\it In particular, $P_g$ is strongly elliptic, formally self-adjoint.}

\vskip 0.18 true cm

We need the method of pseudodifferential operators. If $W$ is an open subset  of ${\Bbb R}^n$, we denote  by $S^m_{1,0}=S^m_{1,0} (W,
{\Bbb R}^n)$ the set of all $p\in C^\infty (W, \mathbb{R}^n)$ such
that for every compact set $O\subset W$ we have
 \begin{eqnarray} \label {-2.3} |D^\beta_x D^\alpha_\xi p(x,\xi)|\le C_{O,\alpha,
 \beta}(1+|\xi|)^{m-|\alpha|}, \quad \; x\in O,\,\, \xi\in {\Bbb R}^n\end{eqnarray}
 for all $\alpha, \beta\in {\Bbb N}^n$, where $|\xi|= \big(\sum_{j=1}^n \xi_j^2\big)^{1/2}$, $D^\alpha =D_1^{\alpha_1}\cdots D_n^{\alpha_n}$, $D_j=\frac{1}{i} \frac{\partial}{\partial x_j}$, ${\mathbb{N}}^n$ is the set of ${\gamma}=(\gamma_1, \cdots, \gamma_n)$ with $\gamma_j=\mbox{integer}\ge 0$, and $|{\gamma}|=\gamma_1+\cdots +\gamma_n$. The  elements  of $S^m_{1,0}$  are  called  symbols (or full symbols) of order $m$.
 It is clear that $S^m_{1,0}$ is a
Fr\'{e}chet space with semi-norms given by the smallest constants
which can be used in (\ref{-2.3}) (i.e.,
\begin{eqnarray*} \|p\|_{O,\alpha, \beta}=
 \,\sup_{x\in O}\big|\left(D_x^\beta D_\xi^\alpha
 p(x, \xi)\right)(1+|\xi|)^{|\alpha|-m}\big|).\end{eqnarray*}
 Let $p(x, \xi)\in S^m_{1,0}$. A pseudodifferential operator
in an open set $W\subset {\Bbb R}^n$ is essentially defined by a
Fourier integral  operator (cf. \cite{KN}, \cite{Ho4}, \cite{Ta2}, \cite{Gr}):
\begin{eqnarray} \label{-2.1}  p(x,D) u(x) = \frac{1}{(2\pi)^{n}} \int_{{\Bbb R}^n}
 e^{i  x\cdot \xi} p(x,\xi) \hat{u} (\xi)d\xi,\end{eqnarray}
and denoted by $OPS^m$.
Here $u\in C_0^\infty (W)$ and $\hat{u} (\xi)=
 \int_{{\Bbb R}^n} e^{-iy\cdot \xi} u(y)dy$ is the Fourier
transform of $u$.
 If there are smooth $p_{m-j} (x, \xi)$, homogeneous in $\xi$ of degree $m-j$ for $|\xi|\ge 1$, that is,
  $p_{m-j} (x, r\xi) =r^{m-j} p_{m-j} (x, \xi)$ for $r, \,|\xi|\ge 1$, and if
 \begin{eqnarray} \label{-2.2} p(x, \xi) \sim \sum_{j\ge 0} p_{m-j} (x, \xi)\end{eqnarray}
 in the sense that
 \begin{eqnarray}  p(x, \xi)-\sum_{j=0}^k p_{m-j} (x, \xi) \in S^{m-k-1}_{1, 0} \end{eqnarray}
 for all $k$, then we say $p(x,\xi) \in S_{cl}^m$, or just $p(x, \xi)\in S^m$. We call
 $p_m(x, \xi)$ the principal symbol of $p(x, D)$. Sometimes we denote by $\iota(p(x,D))$ the (full) symbol of $p(x,D)$.

  Let ${\Omega}$ be a smooth $n$-dimensional Riemannian
 manifold (of class $C^\infty$). We denote by $C^\infty({\Omega})$ and $C_0^\infty(\Omega)$
 the space of all smooth complex-valued functions on $\Omega$  and the subspace of all
 functions with compact support, respectively. Assume that we are given a linear
 operator
 \begin{eqnarray*} P: C^\infty_0(\Omega) \to C^\infty(\Omega).\end{eqnarray*}
 If $G$ is some chart in $\Omega$ (not necessarily connected) and $\kappa: G\to U$ its
diffeomorphism onto an open set $U \subset {\Bbb R}^n$, then let ${\tilde P}$ be defined by the diagram
\begin{eqnarray*}
\begin{CD}
C_0^\infty (G) @> P>>  C^\infty(G) \\
@A\kappa^* AA @AA\kappa^* A \\
 C_0^\infty(U) @> \tilde P>> C^\infty(U)
\end{CD}
\end{eqnarray*}
where
  $\kappa^*$ is the induced transformation
from $C^\infty(U)$ into $C^\infty(G)$, taking a
function $u$ to the function $u\circ \kappa$.
  (note, in the upper row is the operator $r_G\circ P\circ i_G$, where $i_G$ is the natural
embedding $i_G: C_0^\infty (G)\to C^\infty_0 (\Omega)$ and $r_G$ is the natural restriction $r_G: \, C^\infty(\Omega)\to
C^\infty(G)$; for brevity we denote this operator by the same letter $P$ as the
original operator).
 An operator $P: C_0^\infty (\Omega) \to C^\infty(\Omega)$ is called a pseudodifferential
operator on $\Omega$ if for any chart diffeomorphism $\kappa: G\to U$, the
operator $\tilde P$ defined above is a pseudodifferential
operator on $U$. We denote by $OPS^{m}$ the pseudodifferential operator $P$  of order $m$. We also write $OPS^{-\infty}= \bigcap_{m} OPS^m$.

\vskip 0.16 true cm

It is well known (see \cite{Ho3}, \cite{Ho4} or p.$\,$13 of \cite{Ta2}) that if $p_j(x, D)\in OPS^{m_j}$, $j=1,2$, then \begin{eqnarray*} p_1(x, D) p_2(x, D)=q(x,D) \in OPS^{m_1+m_2},\end{eqnarray*}
 and \begin{eqnarray}\label{18-3-27-1} q(x, \xi)= \sum_{\alpha \ge 0} \frac{i^{|\alpha|}}{\alpha!} D_\xi^\alpha p_1 (x, \xi)\, D^\alpha_x p_2(x, \xi).\end{eqnarray}
 An operator $p(x,D)$ is said to be an elliptic pseudodifferential operator of order $m$ if
  for every compact $O\subset \Omega$ there exists a  positive constant $c=c(O)$ such that \begin{eqnarray*} |p(x, \xi)|\ge c|\xi|^m, \,x\in O,\, |\xi|\ge 1.\end{eqnarray*}
   If $q(x, D)\in OPS^{-m}$ is a pseudodifferential operator of order $-m$ such that
\begin{eqnarray*} q(x,D)p(x,D)=I\;\; \mbox{mod}\;\; OPS^{-\infty},\\
 p(x,D)q(x,D)=I\;\; \mbox{mod}\;\; OPS^{-\infty}, \end{eqnarray*}
 then we say that $q(x,D)$ is a (two-sided) parametrix for $p(x,D)$.
 Furthermore, if $P$ is a non-negative elliptic pseudodifferential operator of order $m$, then the spectrum of $P$ lies in a right half-plane and has a finite lower bound $\rho(P) =
\inf\{\mbox{Re}\, \tau\big| \tau\in \sigma(P)\}$, where $\sigma(P)$ denotes the spectrum of $P$. We can modify the principal symbol $h_m(x, \xi)$ for small $\xi$ such that
$h_m(x,\xi)$  has a positive
lower bound throughout and lies in $\{\tau=re^{i\theta} \big| r>0, |\theta|\le \theta_0\}$, where $\theta_0\in (0,\frac{\pi}{2})$.
According to \cite{Gr}, the resolvent $(\tau-P)^{-1}$ exists and is holomorphic in $\tau$ on a neighborhood of a set
\begin{eqnarray*}  W_{r_0,\epsilon} =\{\tau\in {\Bbb C} \big| |\tau|\ge r_0, \mbox{arg}\, \tau \in [\theta_0+\epsilon, 2\pi -\theta_0-\epsilon], \,\mbox{Re}\, \tau \le \rho (P)-\epsilon\}\end{eqnarray*}
(with $\epsilon>0$). There exists a parametrix $Q'_\tau$ on a neighborhood of a possibly larger set
(with $\delta>0,\epsilon>0$)
\begin{eqnarray*} V_{\delta,\epsilon} =\{ \tau \in {\Bbb C} \big| |\tau| \ge \delta \;\;\mbox{or arg}\, \tau \in [\theta_0+\epsilon, 2\pi-\theta_0-\epsilon]\}\end{eqnarray*}
such that this parametrix coincides with $(\tau-P)^{-1}$ on the intersection. Its symbol $q(x,\xi, \tau)$
in local coordinates is holomorphic in $\tau$ there and has the form (cf. Section 3.3 of \cite{Gr})
\begin{eqnarray} \label {-2.6} q(x, \xi,\tau) \sim \sum_{l\ge 0} q_{-m-l} (x, \xi, \tau),\end{eqnarray}
where \begin{eqnarray} \label{8-2.7} & q_{-m}= (\tau- p_m(x, \xi))^{-1}, \quad \;  q_{-m-1} = b_{1,1}(x, \xi) q^2_{-m},\\
& \, \cdots, \,  q_{-m-l} = \sum_{k=1}^{2l} b_{l,k} (x, \xi) q^{k+1}_{-m}, \cdots, l\ge 2  \nonumber\end{eqnarray}
with symbols $b_{l,k}$ independent of $\tau$  and homogeneous of degree $mk-l$ in $\xi$ for $|\xi|\ge1$.
 The semigroup $(e^{-tP})_{t\ge 0}$ can be defined from $P$ by the Cauchy integral formula (see p.$\,$4 of \cite{GG}):
 \begin{eqnarray*} e^{-tP} =\frac{1}{2\pi i} \int_{\mathcal{C}} e^{-t\tau} (\tau-P)^{-1} d\tau, \quad t\ge 0,\end{eqnarray*}
where $\mathcal{C}$ is a suitable curve in the complex plane in the positive direction around the spectrum of $P$.
Inserting (\ref{-2.6})--(\ref{8-2.7}) into above formula, we  get the symbol $\frac{1}{2\pi i} \int_{\mathcal{C}} e^{-t\tau} \big[
\sum_{l\ge 0} q_{-m-l} (x, \xi, \tau)\big]\, d\tau$ of the semigroup $(e^{-tP})_{t\ge 0}$, and furthermore we can  obtain the semigroup $(e^{-tP})_{t\ge 0}$ and its trace for any fixed $t\ge 0$.
\vskip 0.12 true cm

\vskip 1.69 true cm

\section{Full symbol of resolvent operator $(\tau {I}-P_g)^{-1}$}

\vskip 0.45 true cm

Let $(\Omega,g)$ be an $n$-dimensional Riemannian manifold with metric $g=(g_{ij})$.
Note that for $\mathbf{u}=\sum\limits_{j=1}^n u^j \frac{\partial }{\partial x_j}\in T\Omega$, $$\Delta_g u^j= \sum\limits_{m,l=1}^n \Big(g^{ml} \frac{\partial^2 u^j}{\partial x_m\partial x_l} -\sum_{s=1}^n g^{ml} \Gamma^s_{ml} \frac{\partial u^j}{\partial x_s}\Big)$$ and
\begin{eqnarray*}  \mbox{grad}\, \mbox{div}\, \mathbf{u}= \sum\limits_{j,k,m=1}^n \Big(g^{jm}  \big( \frac{\partial^2 u^k}{\partial x_m\partial x_k} +\sum\limits_{l=1}^n \Gamma_{kl}^l \frac{\partial u^k}{\partial x_m} +\sum\limits_{l=1}^n \frac{\partial \Gamma_{kl}^l}{\partial x_m} u^k\big) \Big) \frac{\partial}{\partial x_j}. \end{eqnarray*}
By Lemma 2.1 we can write the Navier-Lam\'{e} operator $P_g$ in $\Omega$  as the form of components relative to coordinates:
 \begin{eqnarray*}  P_g\mathbf{u}=\!\!\!\!\!\!&&\!\!\!\left\{ -\mu\Big( \sum_{m,l=1}^n g^{ml} \frac{\partial^2 }{\partial x_m\partial x_l}\Big){\mathbf{I}}_n -(\mu+\lambda)
\begin{bmatrix} \sum\limits_{m=1}^n g^{1m} \frac{\partial^2}{ \partial x_m\partial x_1} & \cdots & \sum\limits_{m=1}^n g^{1m} \frac{\partial^2}{ \partial x_m\partial x_n}
  \\    \vdots & {} & \vdots  \\
 \sum\limits_{m=1}^n g^{nm} \frac{\partial^2}{ \partial x_m\partial x_1} & \cdots & \sum\limits_{m=1}^n g^{nm} \frac{\partial^2}{ \partial x_m\partial x_n}   \end{bmatrix} \right.\\
&& +\mu \Big( \sum\limits_{m,l,s=1}^n g^{ml} \Gamma_{ml}^s \frac{\partial}{\partial x_s} \Big){\mathbf{I}}_n - \mu  \begin{bmatrix}  \sum\limits_{m,l=1}^n 2g^{ml} \Gamma_{1m}^1 \frac{\partial }{\partial x_l} & \cdots &  \sum\limits_{m,l=1}^n 2g^{ml} \Gamma_{nm}^1 \frac{\partial }{\partial x_l}\\
\vdots & {}& \vdots \\
\sum\limits_{m,l=1}^n 2g^{ml} \Gamma_{1m}^n \frac{\partial }{\partial x_l} & \cdots & \sum\limits_{m,l=1}^n 2g^{ml} \Gamma_{nm}^n \frac{\partial }{\partial x_l}
\end{bmatrix} \\
&& -(\mu+\lambda) \begin{bmatrix}\sum\limits_{m,l=1}^n g^{1m} \Gamma_{1l}^l  \frac{\partial }{\partial x_m} & \cdots & \sum\limits_{m,l=1}^n g^{1m} \Gamma_{nl}^l  \frac{\partial }{\partial x_m} \\       \vdots & {} & \vdots  \\
  \sum\limits_{m,l=1}^n g^{nm} \Gamma_{1l}^l  \frac{\partial }{\partial x_m} & \cdots & \sum\limits_{m,l=1}^n g^{nm} \Gamma_{nl}^l  \frac{\partial }{\partial x_m} \end{bmatrix}
 \\
 &&  -\mu \begin{bmatrix} \sum\limits_{l,m=1}^n g^{ml} \big(  \frac{\partial \Gamma^1_{1l}}{\partial x_m} +  \Gamma_{hl}^1 \Gamma_{1m}^h - \Gamma_{1h}^1 \Gamma_{ml}^h \big) & \cdots &  \sum\limits_{l,m=1}^n g^{ml} \big(  \frac{\partial \Gamma^1_{nl}}{\partial x_m} +  \Gamma_{hl}^1 \Gamma_{nm}^h - \Gamma_{nh}^1 \Gamma_{ml}^h \big)
\\ \vdots & {} & \vdots \\
\sum\limits_{l,m=1}^n g^{ml} \big(  \frac{\partial \Gamma^n_{1l}}{\partial x_m} +  \Gamma_{hl}^n \Gamma_{1m}^h - \Gamma_{1h}^n \Gamma_{ml}^h \big) & \cdots &  \sum\limits_{l,m=1}^n g^{ml} \big(  \frac{\partial \Gamma^n_{nl}}{\partial x_m} +  \Gamma_{hl}^n \Gamma_{nm}^h - \Gamma_{nh}^n \Gamma_{ml}^h \big)
\end{bmatrix} \\
&& \left.- (\mu+\lambda) \begin{bmatrix} \sum_{l,m=1}^n g^{1m} \frac{\partial \Gamma^l_{1l}}{\partial x_m}& \cdots & \sum_{l,m=1}^n g^{1m} \frac{\partial \Gamma^l_{nl}}{\partial x_m}\\
\vdots & {} & \vdots \\
\sum_{l,m=1}^n g^{nm} \frac{\partial \Gamma^l_{1l}}{\partial x_m}& \cdots & \sum_{l,m=1}^ng^{nm} \frac{\partial \Gamma^l_{nl}}{\partial x_m}\end{bmatrix} - \mu \begin{bmatrix}
R^1_1 & \cdots & R^1_n\\
\vdots & {} & \vdots \\
R^n_1 & \cdots & R^n_n\end{bmatrix} \right\}\begin{bmatrix} u^1\\
\vdots\\
u^n\end{bmatrix},   \end{eqnarray*}
where ${\mathbf{I}}_n$ is the $n\times n$ identity matrix. Furthermore, we have
\begin{eqnarray*} P_g\mathbf{u} (x)=\frac{1}{(2\pi)^n} \int_{{\Bbb R}^n}  e^{i  x\cdot \xi}({\mathbf{A}}_g (x, \xi)) \begin{pmatrix} {\hat{u}}_1(\xi)\\
\vdots \\
{\hat{u}}_n (\xi)\end{pmatrix} d\xi,\end{eqnarray*}
where \begin{eqnarray}\label{2020.7.2-1}  && {\mathbf{A}}_g (x, \xi)=
 \mu\Big( \sum_{m,l=1}^n g^{ml} \xi_m\xi_l \Big){\mathbf{I}}_n +(\mu+\lambda)
\begin{bmatrix} \sum\limits_{m=1}^n g^{1m} \xi_m\xi_1 & \cdots & \sum\limits_{m=1}^n g^{1m} \xi_m\xi_n
  \\    \vdots & {} & \vdots  \\
 \sum\limits_{m=1}^n g^{nm} \xi_m\xi_1  & \cdots & \sum\limits_{m=1}^n g^{nm} \xi_m\xi_n  \end{bmatrix}\\
 && \;\; \quad +\text{i}\,\mu \Big( \sum\limits_{m,l,s=1}^n g^{ml} \Gamma_{ml}^s \xi_s \Big){\mathbf{I}}_n - \text{i}\,\mu  \begin{bmatrix}  \sum\limits_{m,l=1}^n 2g^{ml} \Gamma_{1m}^1 \xi_l & \cdots &  \sum\limits_{m,l=1}^n 2g^{ml} \Gamma_{nm}^1 \xi_l\\
\vdots & {}& \vdots \\
\sum\limits_{m,l=1}^n 2g^{ml} \Gamma_{1m}^n \xi_l  & \cdots & \sum\limits_{m,l=1}^n 2g^{ml} \Gamma_{nm}^n \xi_l \end{bmatrix}\nonumber
\\ && \quad \;\; -\text{i}\,(\mu+\lambda) \begin{bmatrix}\sum\limits_{m,l=1}^n g^{1m} \Gamma_{1l}^l  \xi_m & \cdots & \sum\limits_{m,l=1}^n g^{1m} \Gamma_{nl}^l  \xi_m \\       \vdots & {} & \vdots  \\
  \sum\limits_{m,l=1}^n g^{nm} \Gamma_{1l}^l  \xi_m & \cdots & \sum\limits_{m,l=1}^n g^{nm} \Gamma_{nl}^l  \xi_m \end{bmatrix}\nonumber\\
  &&  \quad \;\;-\mu \!\begin{bmatrix} \sum\limits_{l,m=1}^n \!\!g^{ml} \big(  \frac{\partial \Gamma^1_{\!1l}}{\partial x_m}\! +\!  \Gamma_{\!hl}^1 \Gamma_{\!1m}^h \!-\! \Gamma_{\!1h}^1 \Gamma_{\!ml}^h \big) & \cdots &  \sum\limits_{l,m=1}^n \!\!g^{ml} \big(  \frac{\partial \Gamma^1_{\!nl}}{\partial x_m} \!+\!  \Gamma_{\!hl}^1 \Gamma_{nm}^h \!\!-\! \Gamma_{\!nh}^1 \Gamma_{\!ml}^h \big)
\\ \vdots & {} & \vdots \\
\sum\limits_{l,m=1}^n g^{ml}\! \big(  \frac{\partial \Gamma^n_{\!1l}}{\partial x_m}\! +\!  \Gamma_{\!hl}^n \Gamma_{1m}^h \!\!-\! \Gamma_{\!1h}^n \Gamma_{\!ml}^h \big) & \cdots &  \sum\limits_{l,m=1}^n \!\!g^{ml} \big(  \frac{\partial \Gamma^n_{\!nl}}{\partial x_m} \!+\!  \Gamma_{\!hl}^n \Gamma_{\!nm}^h \!\!-\! \Gamma_{\!nh}^n \Gamma_{\!ml}^h \big)\!
\end{bmatrix}\nonumber\\
&& \quad \;\; - (\mu+\lambda) \begin{bmatrix} \sum_{l,m=1}^n g^{1m} \frac{\partial \Gamma^l_{1l}}{\partial x_m}& \cdots & \sum_{l,m=1}^n g^{1m} \frac{\partial \Gamma^l_{nl}}{\partial x_m}\\
\vdots & {} & \vdots \\
\sum_{l,m=1}^n g^{nm} \frac{\partial \Gamma^l_{1l}}{\partial x_m}& \cdots & \sum_{l,m=1}^n g^{nm} \frac{\partial \Gamma^l_{nl}}{\partial x_m}\end{bmatrix} - \mu \begin{bmatrix}
R^1_1 & \cdots & R^1_n\\
\vdots & {} & \vdots \\
R^n_1 & \cdots & R^n_n\end{bmatrix}.\nonumber
  \end{eqnarray}
   For each $\tau\in \mathbb{C}$, we denote \begin{eqnarray} \label{18-3-27-2} \tau {\mathbf{I}}_n- {\mathbf{A}}_g={\mathbf{a}}_2+{\mathbf{a}}_1+{\mathbf{a}}_0,\end{eqnarray} where
   \begin{eqnarray} \label{2020.7.1-1} &&   {\mathbf{a}}_2 (x, \xi):=
 \Big(\!\tau -\mu \sum_{m,l=1}^n g^{ml} \xi_m\xi_l \Big){\mathbf{I}}_n -(\mu+\lambda)
\begin{bmatrix} \sum\limits_{m=1}^n g^{1m} \xi_m\xi_1 & \cdots & \sum\limits_{m=1}^n g^{1m} \xi_m\xi_n
  \\    \vdots & {} & \vdots  \\
 \sum\limits_{m=1}^n g^{nm} \xi_m\xi_1  & \cdots & \sum\limits_{m=1}^n g^{nm} \xi_m\xi_n  \end{bmatrix},\end{eqnarray}
 \begin{eqnarray} \label{2020.7.1-2}
&&  {\mathbf{a}}_1 (x, \xi):=-\text{i}\,\mu \Big( \sum\limits_{m,l,s=1}^n g^{ml} \Gamma_{ml}^s \xi_s \Big){\mathbf{I}}_n + \text{i}\,\mu  \begin{bmatrix}  \sum\limits_{m,l=1}^n 2g^{ml} \Gamma_{1m}^1 \xi_l & \cdots &  \sum\limits_{m,l=1}^n 2g^{ml} \Gamma_{nm}^1 \xi_l\\
\vdots & {}& \vdots \\
\sum\limits_{m,l=1}^n 2g^{ml} \Gamma_{1m}^n \xi_l  & \cdots & \sum\limits_{m,l=1}^n 2g^{ml} \Gamma_{nm}^n \xi_l \end{bmatrix}
\\ && \qquad \qquad\;\;\;\;\,  +\,\text{i}\,(\mu+\lambda) \begin{bmatrix}\sum\limits_{m,l=1}^n g^{1m} \Gamma_{1l}^l  \xi_m & \cdots & \sum\limits_{m,l=1}^n g^{1m} \Gamma_{nl}^l  \xi_m \\       \vdots & {} & \vdots  \\
  \sum\limits_{m,l=1}^n g^{nm} \Gamma_{1l}^l  \xi_m & \cdots & \sum\limits_{m,l=1}^n g^{nm} \Gamma_{nl}^l  \xi_m \end{bmatrix}\nonumber\end{eqnarray}
 \begin{eqnarray}
 && \;\;  {\mathbf{a}}_0 (x, \xi)\!:=\!\mu \!\begin{bmatrix} \sum\limits_{l,m=1}^n \!\!g^{ml} \big(  \frac{\partial \Gamma^1_{\!1l}}{\partial x_m}\! +\!  \Gamma_{\!hl}^1 \Gamma_{\!1m}^h \!-\! \Gamma_{\!1h}^1 \Gamma_{\!ml}^h \big) & \cdots &  \sum\limits_{l,m=1}^n \!\!g^{ml} \big(  \frac{\partial \Gamma^1_{\!nl}}{\partial x_m} \!+\!  \Gamma_{\!hl}^1 \Gamma_{nm}^h \!\!-\! \Gamma_{\!nh}^1 \Gamma_{\!ml}^h \big)
\\ \vdots & {} & \vdots \\
\sum\limits_{\,l,m=1}^n g^{ml}\! \big(  \frac{\partial \Gamma^n_{\!1l}}{\partial x_m}\! +\!  \Gamma_{\!hl}^n \Gamma_{1m}^h \!\!-\! \Gamma_{\!1h}^n \Gamma_{\!ml}^h \big) & \cdots &  \sum\limits_{l,m=1}^n \!\!g^{ml} \big(  \frac{\partial \Gamma^n_{\!nl}}{\partial x_m} \!+\!  \Gamma_{\!hl}^n \Gamma_{\!nm}^h \!\!-\! \Gamma_{\!nh}^n \Gamma_{\!ml}^h \big)\!
\end{bmatrix}\\
&& \quad\qquad\quad \quad + (\mu+\lambda) \begin{bmatrix} \sum_{l,m=1}^n g^{1m} \frac{\partial \Gamma^l_{1l}}{\partial x_m}& \cdots & \sum_{l,m=1}^n g^{1m} \frac{\partial \Gamma^l_{nl}}{\partial x_m}\\
\vdots & {} & \vdots \\
\sum_{l,m=1}^n g^{nm} \frac{\partial \Gamma^l_{1l}}{\partial x_m}& \cdots &\sum_{l,m=1}^n  g^{nm} \frac{\partial \Gamma^l_{nl}}{\partial x_m}\end{bmatrix} + \mu \begin{bmatrix}
R^1_1 & \cdots & R^1_n\\
\vdots & {} & \vdots \\
R^n_1 & \cdots & R^n_n\end{bmatrix}.\nonumber
  \end{eqnarray}

 We will calculate the full symbol of the resolvent operator $(\tau {I} -P_g)^{-1}$. Let $Q$ be a pseudodifferential operator which approximates the resolvent operator
$(\tau {I} -P_g)^{-1}$, i.e.,
\begin{eqnarray*} (\tau I- P_g)Q= I \quad \mbox{mod}\,\, OPS^{-\infty},\\
Q(\tau I- P_g) =I \quad \mbox{mod}\,\, OPS^{-\infty}.\end{eqnarray*}
Let \begin{eqnarray} \label{3--1} {\mathbf{q}}(x, \xi, \tau) \sim {\mathbf{q}}_{-2}(x, \xi, \tau) + {\mathbf{q}}_{-3}(x, \xi, \tau)+\cdots + {\mathbf{q}}_{-2-l}(x, \xi, \tau) +\cdots\end{eqnarray} be the expansion of the full symbol of $Q$.
 Suppose that the complex parameter $\tau$ have homogeneity $2$ (This point of view stems from \cite{See} or \cite{Gil}). Let ${\mathbf{q}}_{-2-l}(x, \xi, \tau)$ be homogeneous of order $-2-l$ in the variables $(\xi, \tau^{1/2})$. This infinite sum defines $\mathbf{q}(x, \xi, \tau)$ asymptotically. Our purpose is to determine $\mathbf{q}(x, \xi, \tau)$ so that
    \begin{eqnarray}\label{3/1}  \iota((\tau I- P_g)Q)\sim {\mathbf{I}}_n,\end{eqnarray}
 where $\iota(T)$ denotes the full symbol of pseudodifferential operator $T$.
 By symbol formula (\ref{18-3-27-1}) of the product of pseudodifferential operators, we can decompose the left-hand side of (\ref{3/1}) into a sum of orders of homogeneity
      \begin{eqnarray}\label {3/2}  \sum_{\alpha\ge 0} \big(\partial_\xi^\alpha (\iota (\tau I- P_g))\big)
            \cdot (D_x^\alpha \mathbf{q})/\alpha !\,\sim {\mathbf{I}}_n, \end{eqnarray}
            where $\partial_\xi^\alpha:=\frac{\partial^{|\alpha|}}{\partial \xi^\alpha}$.
     Noticing that $\iota (\tau I-P_g)=(\tau \mathbf{I}- {\mathbf{A}}_g)= {\mathbf{a}}_2+{\mathbf{a}}_1+{\mathbf{a}}_0$, we find by (\ref{3/2}) that
 \begin{eqnarray*}  \iota ((\tau I- P_g)Q)\sim \sum_{l=0}^\infty \bigg( \sum_{l=j+|\alpha|+2-k} (\partial_\xi^\alpha {\mathbf{a}}_k)\cdot (D_x^\alpha {\mathbf{q}}_{-2-j})/\alpha !\bigg).\end{eqnarray*} The sum is over terms which are homogeneous of order $-l$.
    Thus (\ref{3/1}) leads to the following equations
    \begin{eqnarray}  {\mathbf{I}}_n &=& \sum_{0=j+|\alpha|+2-k} (\partial_\xi^\alpha {\mathbf{a}}_k) (D_x^\alpha {\mathbf{q}}_{-2-j})/\alpha ! = {\mathbf{a}}_2 {\mathbf{q}}_{-2},\nonumber\\
   \label{18-4-2-9} 0&=& \sum_{\underset{l\ge 1}{l=j+|\alpha|+2-k}} (\partial^\alpha_\xi {\mathbf{a}}_k)(D_x^\alpha  {\mathbf{q}}_{-2-j} )/\alpha !\\
    &=&  {\mathbf{a}}_2 {\mathbf{q}}_{-2-l}+  \sum_{\underset {j<l} {l=j+|\alpha|+2-k}} (\partial^\alpha_\xi{\mathbf{a}}_k)(D_x^\alpha
     {\mathbf{q}}_{-2-j})/\alpha !,\quad\quad l\ge 1.\nonumber\end{eqnarray}
     These equations determine the ${\mathbf{q}}_{-2-l}$ inductively. In other words,
     we have \begin{eqnarray}\label{2020.7.3-11} && {\mathbf{q}}_{-2} ={\mathbf{a}}_2^{-1}, \;\;\mbox{and}\;\; {\mathbf{q}}_{-2-l} = -{\mathbf{a}}_2^{-1} \bigg( \sum_{j<l} (\partial^\alpha_\xi  {\mathbf{a}}_k)(D^\alpha_x {\mathbf{q}}_{-2-j})/\alpha !\bigg)\;\;  \mbox{for}\;\; l= j+|\alpha|+2-k\ge 1.\end{eqnarray}

In order to calculate the ${\mathbf{a}}_2^{-1}$, by a direct calculation we find  that
\begin{eqnarray*} &&\begin{bmatrix} \sum\limits_{r=1}^n g^{1r} \xi_r \xi_1 &\cdots &  \sum\limits_{r=1}^n g^{1r} \xi_r \xi_n\\
\vdots& {} &\vdots \\
 \sum\limits_{r=1}^n g^{nr} \xi_r \xi_1 &\cdots &  \sum\limits_{r=1}^n g^{nr} \xi_r \xi_n\end{bmatrix}
  \begin{bmatrix} \sum\limits_{m=1}^n g^{1m} \xi_m \xi_1 &\cdots &  \sum\limits_{=1}^n g^{1m} \xi_m \xi_n\\
\vdots& {} &\vdots \\
 \sum\limits_{m=1}^n g^{nm} \xi_m \xi_1 &\cdots &  \sum\limits_{m=1}^n g^{nm} \xi_m \xi_n\end{bmatrix}\\
&&\qquad \;   = \Big(\sum\limits_{l,m=1}^n g^{lm}\xi_l \xi_m \Big) \begin{bmatrix} \sum\limits_{r=1}^n g^{1r} \xi_r \xi_1 &\cdots &  \sum\limits_{r=1}^n g^{1r} \xi_r \xi_n\\
\vdots& {} &\vdots \\
 \sum\limits_{r=1}^n g^{nr} \xi_r \xi_1 &\cdots &  \sum\limits_{r=1}^n g^{nr} \xi_r \xi_n\end{bmatrix}.\end{eqnarray*}
  Thus, the following two matrices play a key role:
  \begin{eqnarray*} F:=\left\{ {\mathbf{I}}_n, \;\; \begin{bmatrix} \sum\limits_{r=1}^n g^{1r} \xi_r \xi_1 &\cdots &  \sum\limits_{r=1}^n g^{1r} \xi_r \xi_n\\
\vdots& {} &\vdots \\
 \sum\limits_{r=1}^n g^{nr} \xi_r \xi_1 &\cdots &  \sum\limits_{r=1}^n g^{nr} \xi_r \xi_n\end{bmatrix}\right\}.\end{eqnarray*}
The set $F$ of above two matrices can generate a matrix ring $\mathfrak{F}$ according to the usual matrix addition and multiplication of $\mathfrak{F}$
on the ring of functions.
This implies that ${\mathbf{a}}_2^{-1}$ must have the following form:
\begin{eqnarray} \label{2020.7.3-8} {\mathbf{a}}_2^{-1} =  s_1 {\mathbf{I}}_n +s_2\begin{bmatrix} \sum\limits_{r=1}^n g^{1r} \xi_r \xi_1 &\cdots &  \sum\limits_{r=1}^n g^{1r} \xi_r \xi_n\\
\vdots& {} &\vdots \\
 \sum\limits_{r=1}^n g^{nr} \xi_r \xi_1 &\cdots &  \sum\limits_{r=1}^n g^{nr} \xi_r \xi_n\end{bmatrix},\end{eqnarray}
   where $s_1$ and $s_2$ are unknown functions which will be determined later. This key idea is inspired by Galois group theory for solving the polynomial equation (see \cite{Art} or \cite{EHM}).
 Substituting (\ref{2020.7.3-8}) into $ {\mathbf{a}}_2{\mathbf{a}}_2^{-1}={\mathbf{I}}_n$, we have
 \begin{eqnarray*}&& \left\{\Big(\!\tau -\mu \sum_{m,l=1}^n g^{ml} \xi_m\xi_l \Big){\mathbf{I}}_n -(\mu+\lambda)
\begin{bmatrix} \sum\limits_{m=1}^n g^{1m} \xi_m\xi_1 & \cdots & \sum\limits_{m=1}^n g^{1m} \xi_m\xi_n
  \\    \vdots & {} & \vdots  \\
 \sum\limits_{m=1}^n g^{nm} \xi_m\xi_1  & \cdots & \sum\limits_{m=1}^n g^{nm} \xi_m\xi_n  \end{bmatrix}\right\}\Bigg\{  s_1 {\mathbf{I}}_n\Bigg. \\
&&\left.\qquad \quad  +s_2\begin{bmatrix} \sum\limits_{r=1}^n g^{1r} \xi_r \xi_1 &\cdots &  \sum\limits_{r=1}^n g^{1r} \xi_r \xi_n\\
\vdots& {} &\vdots \\
 \sum\limits_{r=1}^n g^{nr} \xi_r \xi_1 &\cdots &  \sum\limits_{r=1}^n g^{nr} \xi_r \xi_n\end{bmatrix}\right\}={\mathbf{I}}_n,\end{eqnarray*}
 i.e., \begin{eqnarray*}&& s_1 \Big(\tau- \mu \sum\limits_{l,m=1}^n g^{lm}\xi_l \xi_m \Big) {\mathbf{I}}_n +\bigg\{s_2 \Big(\tau- \mu \sum\limits_{l,m=1}^n g^{lm}\xi_l \xi_m \Big)- s_1 (\mu+\lambda) \\
 && \qquad \quad -s_2 (\mu+\lambda) \sum\limits_{l,m=1}^n g^{lm} \xi_l \xi_m\bigg\} \begin{bmatrix} \sum\limits_{r=1}^n g^{1r} \xi_r \xi_1 &\cdots &  \sum\limits_{r=1}^n g^{1r} \xi_r \xi_n\\
\vdots& {} &\vdots \\
 \sum\limits_{r=1}^n g^{nr} \xi_r \xi_1 &\cdots &  \sum\limits_{r=1}^n g^{nr} \xi_r \xi_n\end{bmatrix}={\mathbf{I}}_n. \end{eqnarray*}
 Since the set  $F$ is a basis of the matrix ring $\mathfrak{F}$, we get
 \begin{eqnarray*}\left\{ \begin{array}{ll} s_1  \Big(\tau- \mu \sum\limits_{l,m=1}^n g^{lm}\xi_l \xi_m \Big) =1,\\
 s_2 \Big(\tau- \mu \sum\limits_{l,m=1}^n g^{lm}\xi_l \xi_m \Big)-s_1 (\mu+\lambda) -s_2 (\mu+\lambda) \sum\limits_{l,m=1}^n g^{lm}\xi_l \xi_m=0.\end{array}\right.\end{eqnarray*}
 It follows that
 \begin{eqnarray}\left\{ \begin{array}{ll} s_1= \frac{1}{\tau- \mu \sum\limits_{l,m=1}^n g^{lm}\xi_l \xi_m },\\
 s_2= \frac{\mu+\lambda}{ \big(\tau- \mu \sum\limits_{l,m=1}^n g^{lm}\xi_l \xi_m \big)\big(\tau- (2\mu+\lambda)  \sum\limits_{l,m=1}^n g^{lm}\xi_l \xi_m \big)}.\end{array}\right.\end{eqnarray}
Therefore \begin{eqnarray} \label{2020.7.3-10}&& {\mathbf{q}}_{-2} ={\mathbf{a}}_2^{-1}= \frac{1}{\tau- \mu \sum\limits_{l,m=1}^n g^{lm}\xi_l \xi_m }\,{\mathbf{I}}_n\\
 && \qquad \;\; \; \; +\frac{\mu+\lambda}{  \big(\tau- \mu \sum\limits_{l,m=1}^n g^{lm}\xi_l \xi_m \big)\big(\tau- (2\mu+\lambda)  \sum\limits_{l,m=1}^n g^{lm}\xi_l \xi_m \big)}\begin{bmatrix} \sum\limits_{r=1}^n g^{1r} \xi_r \xi_1 &\cdots &  \sum\limits_{r=1}^n g^{1r} \xi_r \xi_n\\
\vdots& {} &\vdots \\
 \sum\limits_{r=1}^n g^{nr} \xi_r \xi_1 &\cdots &  \sum\limits_{r=1}^n g^{nr} \xi_r \xi_n\end{bmatrix}.\nonumber\end{eqnarray}
 Combining this and (\ref{2020.7.3-11}) we can get ${\mathbf{q}}_{-2-l}$ for all $l\ge 1$. For example, we can easily write out the first three terms ${\mathbf{q}}_{-2}$, ${\mathbf{q}}_{-3}$, ${\mathbf{q}}_{-4}$:
\begin{eqnarray}\label{2020.7.3-14} && {\mathbf{q}}_{-2}(x,\xi,\tau)= {\mathbf{a}}_2^{-1},\\
&& {\mathbf{q}}_{-3}(x,\xi,\tau)= -{\mathbf{a}}_2^{-1}\Big({\mathbf{a}}_1  {\mathbf{a}}_2^{-1} - i \sum\limits_{l=1}^n \frac{\partial {\mathbf{a}}_2}{\partial \xi_l}\frac{\partial {\mathbf{a}}_2^{-1}}{\partial x_l} \Big),\\
&& {\mathbf{q}}_{-4}(x,\xi,\tau)=-{\mathbf{a}}_2^{-1}\Big(\sum\limits_{j<2, \, |\alpha|=k-j} \big(\partial^\alpha_\xi \mathbf{a}_k\big) \big( D^\alpha_x \mathbf{q}_{-2-j}\big)/\alpha! \Big).
\end{eqnarray}

 From (\ref{2020.7.3-10}) we immediately have the following:

 \vskip 0.24 true cm

  \noindent{\bf Lemma 3.1.} \ {\it   Let  $Q$ be a pseudodifferential operator satisfy (\ref{3/1}) and let ${\mathbf{q}}_{-2}(x, \xi, \tau)$
   be the principal symbol of $Q$. Then, for any $n\ge 1$,
     \begin{eqnarray}\label{18-4-2-7} && \;\;{\mathbf{q}}_{-2}(x, \xi, \tau)= {\mathbf{a}}_2^{-1},\\
   \label{18-4-2-8}
  && \;\;\mbox{Tr}\,\big({\mathbf{q}}_{-2}(x, \xi, \tau)\big)\!=\!\frac{n}{\big(\tau \!-\!\mu \!\sum_{l,m=1}^n \!g^{lm} \xi_l\xi_m\big)} \!+\! \frac{(\mu+\lambda)\sum_{l,m=1}^n g^{lm} \xi_l\xi_m}{
  \big(\!\tau\! -\!\mu\! \sum_{l,m=1}^n \!g^{lm} \xi_l\xi_m\!\big)
 \big(\!\tau \!-\!(2\mu\!+\!\lambda)  \!\sum_{l,m=1}^n \!g^{lm} \xi_l\xi_m\!\big)},\end{eqnarray}  where ${\mathbf{a}}_2^{-1}$ is given by (\ref{2020.7.3-10}). }

 \vskip 0.20 true cm

\vskip 1.29 true cm

\section{Asymptotic expansion of trace of the integral kernel}

\vskip 0.45 true cm

 \noindent  {\it Proof of Theorem 1.1.} \  From the theory of elliptic operators (see  \cite{Mo1}, \cite{Mo2}, \cite{Mo3}, \cite{Pa}, \cite{Liu-1}, \cite{St}), we see that the Navier-Lam\'{e} operator $-P_g$ can generate  strongly continuous semigroups $(e^{-tP_g^\mp})_{t\ge 0}$ with respect to the Dirichlet and Neumann boundary conditions, respectively, in suitable spaces of vector-valued functions (for example, in $[C_0(\Omega)]^n$ (see \cite{St}) or in $[L^2(\Omega)]^n$ (see \cite{Brow})). Furthermore,
   there exist  matrix-valued functions ${\mathbf{K}}^\mp (t, x, y)$, which are called the integral kernels, such that (see \cite{Brow} or p.$\,$4 of \cite{Frie})
        \begin{eqnarray*}  e^{-tP^\mp_g}{\mathbf{w}}_0(x)=\int_\Omega {\mathbf{K}}^\mp(t, x,y) {\mathbf{w}}_0(y)dy, \quad \,
        {\mathbf{w}}_0\in  [L^2(\Omega)]^n.\end{eqnarray*}

Let $\{{{u}}_k^\mp\}_{k=1}^\infty$ be the orthnormal eigenvectors of the elastic operators $P_g^\mp$ corresponding to the eigenvalues $\{\tau_k^\mp\}_{k=1}^\infty$, then the integral kernels  ${\mathbf{K}}^\mp(t, x, y)=e^{-t P_g^\mp} \delta(x-y)$ are given by \begin{eqnarray} \label{18/12/18} {\mathbf{K}}^\mp(t,x,y) =\sum_{k=1}^\infty e^{-t \tau_k^\mp} {{u}}_k^\mp(x)\otimes {{u}}_k^\mp(y).\end{eqnarray}
This implies that the integrals of the traces of ${\mathbf{K}}^\mp(t,x,y)$ are actually  spectral invariants:
\begin{eqnarray} \label{1-0a-2}\int_{\Omega} \mbox{Tr}\,({\mathbf{K}}^\mp(t,x,x)) dV=\sum_{k=1}^\infty e^{-t \tau_k^\mp}.\end{eqnarray}

 We will combine calculus of symbols (see \cite{See2}) and ``method of images'' to deal with asymptotic expansions for the integrals of traces of integral kernels.
Let $\mathcal{M}=\Omega \cup (\partial \Omega)\cup \Omega^*$ be the (closed) double of $\Omega$, and $\mathcal{P}$ the double to $\mathcal{M}$ of
 the  operator $P_g$ on $\Omega$.

 Let us explain the double Riemannian manifold $\mathcal{M}$ and the differential operator $\mathcal{P}$ more precisely, and introduce how to get them from the given Riemannian manifold $\Omega$ and the Navier-Lam\'{e} operator $P_g$.
The double of $\Omega$ is the manifold $\Omega \cup_{\mbox{Id}} \Omega$, where $\mbox{Id}: \partial \Omega\to \partial \Omega$ is the identity map of $\partial \Omega$; it is obtained from $\Omega \sqcup \Omega$ by identifying each boundary point in one  copy of $\Omega$
with same boundary point in the other. It is a smooth manifold without boundary, and contains two regular domains diffeomorphic to $\Omega$  (see, p.$\,$226 of \cite{Lee}). When considering the double differential system $\mathcal{P}$ crossing the boundary, we make use of the coordinates as follows.  Let $x'=(x_1, \cdots, x_{n-1})$ be any local coordinates for $\partial \Omega$.  For each point $(x',0)\in \partial \Omega$, let $x_{n}$ denote the parameter along the unit-speed geodesic starting at $(x',0)$ with initial direction given by the inward boundary normal to $\partial \Omega$ (Clearly, $x_n$ is the geodesic distance from the point $(x',0)$ to the point $(x',x_n)$). In such coordinates $x_{n}>0$ in
 $\Omega$, and $\partial \Omega$ is locally characterized by $x_{n}=0$ (see, \cite{LU} or \cite{Ta2}).
 Since the Navier-Lam\'{e} operator is a linear differential operator defined on $\Omega$, it can be further denoted as $P_g:=P(g^{\alpha\beta}(x), g^{\alpha n}(x), g^{n\beta}(x), g^{nn}(x), \frac{\partial}{\partial x_1}, \cdots, \frac{\partial}{\partial x_{n-1}},\frac{\partial }{\partial x_n})$, where $1\le \alpha, \beta \le n-1$.
  Let  $\varsigma: (x_1, \cdots, x_{n-1}, x_n)\mapsto (x_1, \cdots, x_{n-1}, -x_n)$ be the reflection with respect to the boundary $\partial \Omega$ in $\mathcal {M}$ (here we always assume $x_n\ge 0$).  Then  we can get the $\Omega^*$ from the given $\Omega$ and $\varsigma$.
    Now, we discuss the change of the metric $g$ from $\Omega$ to $\Omega^*$ by $\varsigma$. Recall that the Riemannian metric $(g_{ij})$ is given in the local coordinates $x_1, \cdots, x_n$, i.e., $g_{ij}(x_1,\cdots, x_n)$. In terms of the new coordinates $z_1,\cdots, z_n$, with  $x_i=x_i(z_1,\cdots, z_n), \,\, i=1,\cdots, n,$  the same metric is given by the functions $\tilde{g}_{ij} =\tilde{g}_{ij} (z_1, \cdots, z_n)$, where
\begin{eqnarray} \label{2022.10-2} \tilde{g}_{ij}= \frac{\partial x_k}{\partial z_i} g_{kl} \frac{\partial x_l}{\partial z_j}.\end{eqnarray}
  If $\varsigma$ is a coordinate change in a neighborhood intersecting with $\partial \Omega$
  \begin{eqnarray}\label{2022.10.28-1} \left\{ \begin{array}{ll} x_1= z_1, \\
                   \cdots \cdots\\
                   x_{n-1}=z_{n-1},\\
                   x_n=-z_n,\end{array}\right.\end{eqnarray}
                                       then its Jacobian matrix is
\begin{eqnarray} \label{2022.11.8-1} J:=\begin{pmatrix} \frac{\partial x_1}{\partial z_1}& \cdots &   \frac{\partial x_1}{\partial z_{n-1}} & \frac{\partial x_1}{\partial z_{n}}\\
   \vdots & \ddots & \vdots & \vdots\\
  \frac{\partial x_{n-1}}{\partial z_{1}} & \cdots & \frac{\partial x_{n-1}}{\partial z_{n-1}} &  \frac{\partial x_{n-1}}{\partial z_{n}}\\
   \frac{\partial x_n}{\partial z_{1}}& \cdots & \frac{\partial x_n}{\partial z_{n-1}}&\frac{\partial x_n}{\partial z_{n}}\end{pmatrix} =   \begin{pmatrix} 1& \cdots &   0 & 0\\
   \vdots & \ddots & \vdots & \vdots\\
  0 & \cdots & 1 &  0\\
   0& \cdots & 0&-1\end{pmatrix}.\end{eqnarray}
Using this and (\ref{2022.10-2}), we immediately obtain the corresponding metric on the $\Omega^*$:
  (see \cite{Liu1}, \cite{LiuTan} or p.\,10169, p.\,10183 and p.\,10187 of \cite{Liu-21} )
\begin{eqnarray} \label{2021.2.6-3}  g_{jk} (\overset{*}{x})\!\!\!&\!=\!&\!\!\!- g_{jk} (x) \quad \, \mbox{for}\;\;
  j<k=n \;\;\mbox{or}\;\; k<j=n,\\ g(\overset{*}{x}) \!\!\!&\!=\!&\!\!\! g_{jk} (x)\;\; \;\;\mbox{for}\;\; j,k<n \;\;\mbox{or}\;\; j=k=n,\\
  \label{2021.2.6-4}  g_{jk}(x)\!\!\!&\!=\!&\!\!\! 0 \;\; \;\mbox{for}\;\; j<k=n \;\;\mbox{or}\;\; k<j=n \;\;\mbox{on}\;\; \partial \Omega,\end{eqnarray} where $x_n(\overset{*}{x})= -x_n (x)$.
  We denote such a new (isometric) metric on $\Omega^*$ as $g^*$.
  It is easy to verify that \begin{eqnarray*} \begin{bmatrix} g_{11}(x) & \cdots & g_{1,n-1}(x)& -g_{1n}(x)\\
  \vdots & \ddots & \vdots & \vdots\\
  g_{n-1,1}(x) & \cdots & g_{n-1,n-1}(x) &- g_{n-1,n}(x)\\
 - g_{n1}(x) & \cdots & -g_{n,n-1}(x) & g_{nn}(x)\end{bmatrix}^{-1}=\begin{bmatrix} g^{11}(x) & \cdots & g^{1,n-1}(x)& -g^{1n}(x)\\
  \vdots & \ddots & \vdots & \vdots\\
  g^{n-1,1}(x) & \cdots & g^{n-1,n-1}(x) &- g^{n-1,n}(x)\\
 - g^{n1}(x) & \cdots & -g^{n,n-1}(x) & g^{nn}(x)\end{bmatrix}, \end{eqnarray*}
  where $[g^{jk}(x)]_{n\times n}$ is the inverse of $[g_{jk}(x)]_{n\times n}$.
    In addition, by this reflection $\varsigma$,
the differential operators $\frac{\partial }{\partial x_1}$, $\cdots$, $\frac{\partial}{\partial x_{n-1}}$, $\frac{\partial }{\partial x_n}$ (defined on  $\Omega$) are changed to $\frac{\partial }{\partial x_1}$, $\cdots$, $\frac{\partial}{\partial x_{n-1}}$, $-\frac{\partial }{\partial x_n}$ (defined on  $\Omega^*$), respectively.
 We define \begin{eqnarray} \label{2022.10.18-2} \mathcal{P}=\left\{\begin{array}{ll} \! P_g \;\;\; \;\;\;\,\mbox{on} \;\, \Omega\\
 \!P^* \;\;\;\; \;\mbox{on} \;\, \Omega^*, \end{array} \right.\end{eqnarray}
where \begin{eqnarray} \label{2022.10.6-8} P^*:=P\Big({g}^{\alpha\beta}(\overset {*}{x}), - {g}^{\alpha n}(\overset{*}{x}),- {g}^{n\beta}(\overset{*}{x}), {g}^{nn}(\overset{*}{x}), \frac{\partial }{\partial x_1}, \cdots, \frac{\partial }{\partial x_{n-1}},- \frac{\partial }{\partial x_n}\Big), \end{eqnarray}
and $\overset{*}{x}=(x',-x_n)\in \Omega^*$.
  Clearly, the differential operator $P^*$ is obtained by $P_g$ and the reflection $\varsigma$, that is, $P^*$ is got if we replace $g^{\alpha\beta}(x)$, $g^{\alpha n}(x)$, $g^{n\beta} (x)$, $g^{nn}(x)$, $\frac{\partial }{\partial x_n}$ by ${g}^{\alpha\beta}(\overset{*}{x})$, $-{g}^{\alpha n}(\overset{*}{x})$, $-{g}^{n\beta} (\overset{*}{x})$, ${g}^{nn}(\overset{*}{x})$, $-\frac{\partial }{\partial x_n}$ in $P_g=P\Big( g^{\alpha\beta} (x)$, $g^{\alpha n}(x)$, $g^{n\beta}(x)$, $g^{nn}(x)$, $\frac{\partial }{\partial x_1}$, $\cdots$, $\frac{\partial }{\partial x_{n-1}}$, $\frac{\partial }{\partial x_n}\Big)$, respectively.  Note that  $g^{\alpha\beta} (\overset{*}{x})=
g^{\alpha\beta} (x)$, $-g^{\alpha n}(\overset{*}{x})= g^{\alpha n}(x)$, $-g^{n\beta}(\overset{*}{x})=g^{n\beta}(x)$ and $g^{nn}(\overset{*}{x})=g^{nn} (x)$.
In view of the metric matrices $g$ and $g^*$ have the same order principal minor determinants, we see that $\mathcal{P}$ is still a linear elliptic differential operator on $\mathcal{M}$.

 Let $\mathbf{K}(t,x,y)$ be the fundamental solution of the parabolic  system
 \begin{eqnarray*} \left\{ \begin{array}{ll} \frac{\partial \mathbf{u}}{\partial t} + \mathcal{P}\mathbf{u}=0 \;\; &\mbox{in}\;\, (0,+\infty)\times \mathcal{M},\\
  \mathbf{u}=\boldsymbol{\phi} \;\; &\mbox{on}\;\; \{0\}\times \mathcal{M}.\end{array}\right.\end{eqnarray*}
  That is, for any $t\ge 0$ and $x,y\in \mathcal{M}$,
\begin{eqnarray}\label{2020.10.28-10}\left\{\begin{array}{ll}
    \frac{\partial \mathbf{K}(t,x,y)}{\partial t} + \mathcal{P}\mathbf{K}(t,x,y)=0 \;\;  &\mbox{for}\;\, t>0, \, x, y\in \mathcal{M},\\
      \mathbf{K}(0,x,y)=\boldsymbol{\delta}(x-y) \;\; &\mbox{for}\;\;  x,y \in \mathcal{M}.\end{array}\right.\end{eqnarray}
   Here the operator $\mathcal{P}$ is acted in the third argument $y$ of $\mathbf{K}(t,x,y)$.

 Clearly, the coefficients occurring in $\mathcal{P}$ jump as $x$ crosses the $\partial \Omega$ (since the extended metric $g$ is $C^0$-smooth on whole $\mathcal{M}$ and $C^\infty$-smooth in $\mathcal{M}\setminus \partial \Omega$), but $\frac{\partial \mathbf{u}}{\partial t}+\mathcal{P} \mathbf{u}=0$ with $\mathbf{u}(0,x)=\boldsymbol{\phi}(x)$ still has a nice fundamental solution $\mathbf{K}$ of class $C^1((0,+\infty)\times \mathcal{M}\times  \mathcal{M}) \cap C^\infty((0,+\infty)\times (\mathcal{M} \setminus \partial \Omega) \times (\mathcal{M}\setminus \partial \Omega))$, approximable even on $\partial \Omega$  by Levi's sum (see \cite{Liu-21}, or another proof below). Now, let us restrict $x,y\in \Omega$.
  It can be verified that
$\mathbf{K}^{-}(t,x,y):=\mathbf{K}(t,x,y)- \mathbf{K}(t,x,\overset{*}{y})$ and $\mathbf{K}^{+}:=\mathbf{K}(t,x,y)+ \mathbf{K}(t,x, \overset{*}{y})$   are the Green functions of \begin{eqnarray*}\left\{ \begin{array}{ll}\frac{\partial \mathbf{u}}{\partial t} +{P}_g\mathbf{u}=0\;\;&\mbox{in}\;\, (0,+\infty) \times \Omega,\\
\mathbf{u}=\boldsymbol{\phi} \;\; &\mbox{on}\;\; \{0\}\times \Omega\end{array}\right.\end{eqnarray*} with zero Dirichlet and Neumann boundary conditions, respectively,
where $y=(y',y_n)$, $y_n\ge 0$, and  $\overset{*}{y}:=\varsigma(y',y_n)=(y^{\prime},-y_n)$.  In other words,
 \begin{eqnarray*}\left\{    \begin{array}{ll} \frac{\partial \mathbf{K}^-(t,x,y)}{\partial t} + {P}_g\mathbf{K}^-(t,x,y)=0,\;\;\; t>0,\, x,\, y\in\Omega,\\
 \mathbf{K}^- (t, x,y)=0, \;\;\;\;  t>0,\; x\in \Omega, \,\;  y\in \partial \Omega,\\
 \mathbf{K}^-(0, x,y)=\boldsymbol{\delta}(x-y), \;\;\;  x,y\in \Omega\end{array} \right. \end{eqnarray*}
and
 \begin{eqnarray*}\left\{  \begin{array}{ll}   \frac{\partial \mathbf{K}^+(t,x,y)}{\partial t} + {P}_g\mathbf{K}^+(t,x,y)=0, \,\;\;\; t>0, \;x, \,y\in \Omega,\\
 \frac{\partial \mathbf{K}^+ (t, x,y)}{\partial \boldsymbol{\nu}}=0, \;\; \,\,t>0,\;  x\in \Omega, \;\;  y\in \partial \Omega,\\
 \mathbf{K}^+(0, x,y)=\boldsymbol{\delta}(x-y),\;\;\;  x,y\in \Omega,\end{array} \right. \end{eqnarray*}
 where $\frac{\partial \mathbf{K}^+}{\partial \boldsymbol{\nu}}:= \mu \big(\nabla \mathbf{K}^+ +(\nabla \mathbf{K}^+)^T\big)\boldsymbol{\nu} +\lambda (\mbox{div}\, \mathbf{K}^+)\boldsymbol{\nu}$ on $\partial \Omega$.
  In fact, for any $t>0$, $x,y \in \Omega$, we have
  $P_g\mathbf{K}(t,x, y)=\mathcal{P}\mathbf{K}(t,x, y)$,  so that
  \begin{eqnarray}\label{2022.10.29-8}\left\{ \begin{array}{ll}
  \Big( \frac{\partial}{\partial t}+{P}_g\Big) \mathbf{K} (t,x, y)=
  \Big( \frac{\partial}{\partial t}+\mathcal{P}\Big) \mathbf{K} (t,x, y)=0,\\
  \mathbf{K}(0, x,y)=\boldsymbol{\delta}(x-y) \end{array}\right.\end{eqnarray}
  by (\ref{2020.10.28-10}).
Noting that the  Jacobian matrix of the reflection $\varsigma$ is $J$ (see (\ref{2022.11.8-1})), it follows from chain rule that for any fixed $t>0$ and $x\in \Omega$, and any $y=(y',y_n)\in \Omega$,
\begin{align*} &\left[{P}_g ( \mathbf{K}(t,x,\overset{*}{y}))\right]\bigg|_{\text{\normalsize evaluated at the point $y$}}\\
&=\left[{P}_g ( \mathbf{K}(t,x,\varsigma(y',y_n)))\right]\bigg|_{\text{\normalsize evaluated at the point $(y',y_n)$}}
 \\
&= \left[{P}_g (\mathbf{K}(t, x, (y',-y_n))\big)\right]\bigg|_{\text{\normalsize evaluated at the point $(y',y_n)$}}\\
& \! =\!
 \begin{small}{\left\{\!\left[P \big(g^{\alpha\beta}(y), g^{\alpha n}(y), g^{n\beta} (x), g^{nn}(x), \frac{\partial}{\partial y_1}, \cdots, \frac{\partial}{\partial y_{n-1}}, \frac{\partial}{\partial y_n}\big)\right]\! \! \mathbf{K} (t, x, (y'\!,-y_n))\!\right\}\!\Bigg|_{\text{\normalsize evaluated at $\!(y'\!,y_n)$}}}\end{small}\\
\!& \! =
 \begin{small}{\left\{\!\left[P \big(g^{\alpha\beta}(\overset{*}{y}), -g^{\alpha n}(\overset{*}{y}), -g^{n\beta}(\overset{*}{y}), g^{nn}(\overset{*}{y}), \frac{\partial}{\partial y_1}, \cdots, \frac{\partial}{\partial y_{n-1}}, -\frac{\partial}{\partial y_n}\big)\right]\!  \mathbf{K}(t,x,y)\!\right\}\Bigg|_{\text{\normalsize evaluated  at $\overset{*}{y}=(y',-y_n)$}}}\end{small} \\
&= P^* ( \mathbf{K}(t, x, \overset{*}{y}))\Big|_{\text{\normalsize evaluated at the point $\overset{*}{y}=(y',-y_n)$}}.  \end{align*}
 That is, the action of $P_g$ to $\mathbf{K}(t,x,\overset{*}{y})$ at the point ${y}=(y',y_n)$ is just the action of $P^*$
 to $\mathbf{K}(t,x,\overset{*}{y})$ at the point $\overset{*}{y}=(y',-y_n)$. Because of  $\varsigma(y',y_n)=(y',-y_n)\in \Omega^*$, we see $$P^*( \mathbf{K}(t,x,\overset{*}{y}))\big|_{\text{\normalsize evaluated at the point $\overset{*}{y}=(y',-y_n)$}}=\mathcal{P}(\mathbf{K} (t,x,\overset{*}{y}))\big|_{\text{\normalsize evaluated at the point $\overset{*}{y}=(y',-y_n)$}}.$$
For any $t>0$, $x\in \Omega$ and $(y',-y_n)\in \Omega^*$, we have $$(\frac{\partial}{\partial t}+ {\mathcal{P}}) (\mathbf{K} (t, x,  (y',-y_n)))=0.$$
In addition, $\mathbf{K} (t, x,  (y',-y_n))= \mathbf{K} (t, x,  \varsigma(y))$ for any $t>0$, $x,y\in \Omega$. By virtue of $x\ne (y',-y_n)$,
 this leads to $\mathbf{K} (0, x,  (y',-y_n))=0$ and
 \begin{eqnarray*} \label{2020.10.29-1} (\frac{\partial}{\partial t}+ {P}_g) \big(\mathbf{K} (t, x,  (y',-y_n))\big)=0\;\; \mbox{for any}\;\,t>0, x\in \Omega \;\,\mbox{and}\;\,(y',-y_n)\in \Omega^*,\end{eqnarray*}
 i.e.,
  \begin{eqnarray} \label{2020.10.29-2}\left\{\! \begin{array}{ll}  (\frac{\partial}{\partial t}+ {{P}_g}) \mathbf{K} (t, x,  \overset{*}{y})=0\;\; \mbox{for any}\;\,t>0, x\in \Omega \;\,\mbox{and}\;\,\overset{*}{y}\in \Omega^*,\\
   \mathbf{K} (0, x,  \overset{*}{y})=0 \;\;\, \mbox{for any}\,\; x,y\in \Omega.\end{array} \right.\end{eqnarray}
 Combining (\ref{2022.10.29-8}) and (\ref{2020.10.29-1}), we obtain that
 \begin{eqnarray} \label{2020.10.29-3} \left\{ \begin{array}{ll} (\frac{\partial}{\partial t}+ {{P}_g}) \Big(\mathbf{K} (t, x, {y})-\mathbf{K} (t, x,  \overset{*}{y})\Big)=0\;\;\mbox{for any}\;\,t>0, \; x,y\in \Omega,\\
 \mathbf{K} (0, x, {y})-\mathbf{K} (0, x,  \overset{*}{y})=\boldsymbol{\delta} (x-y)\;\;\mbox{for any}\;\, x,y\in \Omega.\end{array}\right.\end{eqnarray}
  $\mathbf{K}(t,x, y)$ is $C^1$-smooth with respect to $y$ in $\mathcal{M}$ for any fixed $t>0$ and $x\in \Omega$, so does it on the hypersurface $\partial \Omega$. Therefore, we get that  $\mathbf{K}^-(t,x,y)$ (respectively $\mathbf{K}^+(t,x,y)$) is the Green function in $\Omega$ with the Dirichlet (respectively, Neumann) boundary condition on $\partial \Omega$.

\vskip 0.12 true cm

  To show $C^1$-regularity of the fundamental solution $\mathbf{K}(t,x,y)$, it suffices to prove  $C^{1,1}_{loc}$-regularity for the solution $\mathbf{u}$ of the  elliptic system $\mathcal{P}\mathbf{u}=\mathbf{f}$ in $\mathcal{M}$. This immediately follows from Xiong's result \cite{Xi-11} of $C^{1,1}_{loc}$-regularity for solution of elliptic system
   \begin{eqnarray*}\left\{\!\begin{array}{ll} L\mathbf{u}=\mathbf{f} \;\;\mbox{in}\;\, U,\\
   \mathbf{u}=\boldsymbol{\phi} \;\, \,\mbox{on} \,\; \partial U\end{array}\right.\end{eqnarray*}
   with piecewise uniformly H\"{o}lder continuous coefficients and $\mathbf{f}$ on both sides of a general $(n-1)$-dimensional embedded $C^{1,\alpha}$ hypersurface $S$ (the coefficients might be discontinuous cross this hypersurface), where $L$ is a general elliptic operator of second order, $U$ is a bounded domain in $\mathbb{R}^n$ with smooth boundary, $\boldsymbol{\phi}\in [C(\partial U)]^n$ and $S\cap U\ne \emptyset$. This result (see, Theorem 1.2 of \cite{Xi-11}) can also be applied to our case for $C^{1,1}_{loc}$-regularity of the fundamental solution in an $n$-dimensional Riemannian manifold (see \cite{Liu-21}). When the coefficients of an elliptic (or parabolic) system are piecewise smooth between an $(n-1)$-dimensional hyperplane (might be discontinuous cross this hyperplane), the corresponding $C^{1+\alpha/2,\alpha}$-regularity for solution of an elliptic system was obtained by Dong (see Remark 5 of p.$\,$141 in \cite{Do-12}).

\vskip 0.20 true cm

Therefore, the integral kernels $\mathbf{K}^{\mp}(t,x,y)$
 of $\frac{\partial \mathbf{u}}{\partial t}+P_g^{\mp}\mathbf{u}=0$ can be expressed on $(0,\infty)\times \Omega\times \Omega$ as
 \begin{eqnarray} \label{c4-23} {\mathbf{K}}^{\mp} (t,x,y) =\mathbf{K}(t,x,y)\mp \mathbf{K}(t,x,\overset{\ast} {y}),\end{eqnarray}
 $\overset{*} {y}$ being the double of $y\in \Omega$ (see, p.$\,$53 of \cite{MS}).
 Since the strongly continuous semigroup $(e^{-t\mathcal{P}})_{t\ge 0}$ can also be represented as  \begin{eqnarray*} e^{-t\mathcal{P}} =\frac{1}{2\pi i} \int_{\mathcal{C}} e^{-t\tau} (\tau I- \mathcal{P})^{-1} d\tau,\end{eqnarray*}
where $\mathcal{C}$ is a suitable curve in the complex plane in the positive direction around the spectrum of $\mathcal{P}$ (i.e., a contour around the positive real axis). It follows that
 \begin{eqnarray}\quad \;\;\;\;\, {\mathbf{K}} (t,x,y) = e^{-t\mathcal{P}}\delta(x-y) = \frac{1}{(2\pi)^n} \int_{{\Bbb R}^n} e^{i(x-y)\cdot\xi} \bigg(\!\frac{1}{2\pi i} \! \int_{\mathcal{C}} e^{-t\tau}\; \iota \big((\tau {I} -\mathcal{P})^{\!-1}\big) d\tau\!\bigg) d\xi, \;\; \forall x,y\in \mathcal{M}.\end{eqnarray}

 We claim that \begin{eqnarray} \label{2022.11.11-1} \frac{1}{2\pi i} \int_{\mathcal{C}} (\tau I- \mathcal{P} )^{-1} e^{-t\tau}\,\delta(x-y)\,d\tau  = \frac{1}{2\pi i} \int_{\mathcal{C}} \Big( \int_{\mathbb{R}^n} e^{i(x-y)\cdot \xi} \sum_{j\le -2} \mathbf{q}_j ( x,\xi, \tau) \,d\xi \Big) e^{-t \tau} d\tau.\end{eqnarray}
 In fact, for any smooth vector-valued function $\boldsymbol{\phi}$ with compact  we have \begin{eqnarray*} \big(e^{-t\mathcal{P}} \boldsymbol{\phi} \big)(x) \!\!&\!\!=\!\!&\! \! \Big( \frac{1}{2\pi i} \int_{\mathcal{C}} e^{-t\tau} (\tau I- \mathcal{P} )^{-1}  d\tau \Big) \boldsymbol{\phi}(x)\\
 \!&=\!& \frac{1}{2\pi i} \int_{\mathcal{C}} e^{-t\tau} \Big( \int_{\mathbb{R}^n} e^{i x\cdot \xi} \sum_{j\le -2} \mathbf{q}_j (x, \xi, \tau) \hat{\boldsymbol{\phi}} (\xi) \,d\xi \Big) d\tau.\end{eqnarray*}
On the one hand, from the left-hand side of (\ref{2022.11.11-1}), we get \begin{eqnarray} \label{2022.11.11-2}
&& \int\Big[ \Big( \frac{1}{2\pi i} \int_{\mathcal{C}} (\tau I- \mathcal{P} )^{-1}  e^{-t\tau} d\tau\Big) (\delta (x-y))\Big]\boldsymbol{\phi}(y)dy \\
&&   \;\quad\quad \quad =  \Big( \frac{1}{2\pi i} \int_{\mathcal{C}} (\tau- \mathcal{P})^{-1} e^{-t\tau} d\tau\Big) \boldsymbol{\phi}(x)
= e^{-t\mathcal{P} } \boldsymbol{\phi} (x).\nonumber\end{eqnarray}
On the other hand, from the right-hand side of (\ref{2022.11.11-1}) we obtain
\begin{eqnarray} \label{2022.11.11-3} && \int \Big[ \frac{1}{2\pi i} \int_C \Big(\int_{\mathbb{R}^n} e^{i(x-y)\cdot \xi} \sum_{j\le -2} \mathbf{q}_j ( x, \xi,\tau) d\xi\Big) e^{-t\tau} d\tau \Big] \boldsymbol{\phi}(y)dy \\
&& \quad \quad \quad =\frac{1}{2\pi i} \int_C \Big(\int_{\mathbb{R}^n} e^{ix\cdot \xi} \sum_{j\le -2} \mathbf{q}_j ( x, \xi,\tau) d\xi\Big) e^{-t\tau} d\tau \int e^{-y\cdot \xi} \boldsymbol{\phi}(y) dy \nonumber\\
  && \quad \quad \quad = \frac{1}{2\pi i} \int_C \Big(\int_{\mathbb{R}^n} e^{ix\cdot \xi} \sum_{j\le -2} \mathbf{q}_j (x, \xi,\tau) \hat {\boldsymbol{\phi}} (\xi) d\xi\Big) e^{-t\tau} d\tau
= e^{-t\mathcal{P}} \boldsymbol{\phi} (x).\nonumber\end{eqnarray}
 Thus, the desired identity (\ref{2022.11.11-1}) is asserted by (\ref{2022.11.11-2}) and (\ref{2022.11.11-3}).

  In particular, for every $x\in \Omega$, \begin{eqnarray}\label{2020.7.5-1}&&  {\mathbf{K}} (t,x,x) = e^{-t\mathcal{P}}\delta(x-x) = \frac{1}{(2\pi)^n} \int_{{\Bbb R}^n}\bigg(\frac{1}{2\pi i} \int_{\mathcal{C}} e^{-t\tau}\; \iota \big((\tau {I} -\mathcal{P})^{-1}\big) d\tau\bigg) d\xi\\
    && \qquad\,\; \qquad = \frac{1}{(2\pi)^n} \int_{{\Bbb R}^n}\bigg(\frac{1}{2\pi i} \int_{\mathcal{C}} e^{-t\tau}\; \iota \big((\tau {I} -{P_g})^{-1}\big) d\tau\bigg) d\xi\nonumber\\
    && \qquad\,\;\qquad = \frac{1}{(2\pi)^n} \int_{{\mathbb{R}}^n} \Big( \frac{1}{2\pi i} \int_{\mathcal{C}} e^{-t\tau} \sum_{l\ge 0} q_{-2-l} (x, \xi, \tau) \,d\tau \Big) d\xi,\nonumber\\
\label{2020.7.5-2} &&  {\mathbf{K}} (t,x,\overset{*}{x}) = e^{-t\mathcal{P}}\delta(x-\overset{*}{x}) = \frac{1}{(2\pi)^n} \int_{{\Bbb R}^n} e^{i(x-\overset{*}{x})\cdot\xi} \bigg(\frac{1}{2\pi i} \int_{\mathcal{C}} e^{-t\tau}\; \iota \big((\tau {I} -\mathcal{P})^{-1}\big) d\tau\bigg) d\xi\\
&& \qquad \qquad\;\;\; = \frac{1}{(2\pi)^n} \int_{{\Bbb R}^n} e^{i(x-\overset{*}{x})\cdot\xi} \bigg(\frac{1}{2\pi i} \int_{\mathcal{C}} e^{-t\tau}\; \iota \big((\tau {I} -{\mathcal{P}})^{-1}\big) d\tau\bigg) d\xi\nonumber\\
 && \qquad\;\;\;\qquad = \frac{1}{(2\pi)^n} \int_{{\mathbb{R}}^n} e^{i(x-\overset{*}{x})\cdot \xi} \Big( \frac{1}{2\pi i} \int_{\mathcal{C}} e^{-t\tau} \sum_{l\ge 0} q_{-2-l} (x, \xi, \tau) \,d\tau \Big) d\xi,\nonumber
 \end{eqnarray}
 where $\sum_{l\ge 0} {\mathbf{q}}_{-2-l} (x,\xi,\tau) $ is the full symbol of $(\tau I-P_g)^{-1}$.

 Firstly, from the discussion of previous section, we know that
   \begin{eqnarray}\label{2022.11.9-1}
   && \;\;\;\;{\mathbf{q}}_{-2} (x,\xi,\tau) =\frac{1}{\tau- \mu \sum\limits_{l,m=1}^n g^{lm}\xi_l \xi_m }\,{\mathbf{I}}_n\\
 && \qquad \;\;\;\; \; \; +\frac{\mu+\lambda}{  \big(\tau- \mu \sum\limits_{l,m=1}^n g^{lm}\xi_l \xi_m \big)\big(\tau- (2\mu+\lambda)  \sum\limits_{l,m=1}^n g^{lm}\xi_l \xi_m \big)}\begin{bmatrix} \sum\limits_{r=1}^n g^{1r} \xi_r \xi_1 &\cdots &  \sum\limits_{r=1}^n g^{1r} \xi_r \xi_n\\
\vdots& {} &\vdots \\
 \sum\limits_{r=1}^n g^{nr} \xi_r \xi_1 &\cdots &  \sum\limits_{r=1}^n g^{nr} \xi_r \xi_n\end{bmatrix}\nonumber\end{eqnarray}
 and   \begin{eqnarray}\label{2020.6.6-1} \;\;\quad \quad \,\;\mbox{Tr}\, \big({\mathbf{q}}_{-2} (x,\xi,\tau)\big)=
   \frac{n}{\big(\tau \!-\!\mu \!\sum_{l,m=1}^n \!g^{lm} \xi_l\xi_m\big)} \!+\! \frac{(\mu+\lambda)\sum_{l,m=1}^n g^{lm} \xi_l\xi_m}{
  \big(\!\tau\! -\!\mu\! \sum_{l,m=1}^n \!g^{lm} \xi_l\xi_m\!\big)
 \big(\!\tau \!-\!(2\mu\!+\!\lambda)  \!\sum_{l,m=1}^n \!g^{lm} \xi_l\xi_m\!\big)}
.\end{eqnarray}
       For each $x\in \Omega$, we use a geodesic normal coordinate system centered at this $x$. It follows from \S11 of Chap.1 in \cite{Ta1}
        that in such a coordinate system, $g_{jk}(x)=\delta_{jk}$ and $\Gamma_{jk}^l(x)=0$. Then (\ref{2020.6.6-1}) reduces to
        \begin{eqnarray}\label{2020.7.6-3} \quad \quad \,\;\mbox{Tr} \,\big({\mathbf{q}}_{-2} (x, \xi,\tau)\big)=
     \frac{n}{(\tau -\mu |\xi|^2)}+
      \frac{(\mu+\lambda)|\xi|^2}{(\tau -\mu |\xi|^2) (\tau-(2\mu+\lambda)|\xi|^2)},\end{eqnarray}
  where $|\xi|=\sqrt{\sum_{k=1}^n \xi^2_k}$ for any $\xi\in {\mathbb{R}}^n$.
   By applying the residue theorem (see, for example, Chap.$\,$4, \S5 in \cite{Ahl}) we get
 \begin{eqnarray} \label{3.10} && \frac{1}{2\pi i} \int_{\mathcal{C}} e^{-t\tau} \bigg(\frac{n}{(\tau -\mu |\xi|^2)}+
      \frac{(\mu+\lambda)|\xi|^2}{(\tau -\mu |\xi|^2) (\tau-(2\mu+\lambda)|\xi|^2)} \bigg) d\tau= (n-1) e^{-t\mu|\xi|^2} +  e^{-t(2\mu+\lambda)|\xi|^2}.\end{eqnarray}
       It follows that
      \begin{eqnarray}\label{2020.7.10-1} \frac{1}{(2\pi)^n}\!\!\!\!&\!\!\!&\!\!\!\!\! \!\!\! \!\! \int_{{\mathbb{R}}^n}\Big( \frac{1}{2\pi i} \int_{\mathcal{C}} e^{-t\tau}\, \mbox{Tr}\,({\mathbf{q}}_{-2} (x, \xi,\tau) ) d\tau \Big) d\xi \\
      \!\!\! &=\!\!\!&
    \frac{1}{(2\pi)^n} \int_{{\Bbb R}^n} \bigg((n-1) e^{-t\mu|\xi|^2} +  e^{-t(2\mu+\lambda)|\xi|^2}\bigg)
         d\xi \nonumber\\
               \!\!\! &=\!\!\!&\frac{n-1}{(4\pi \mu t)^{n/2}} +  \frac{1}{(4\pi (2\mu+\lambda) t)^{n/2}},
                  \nonumber\end{eqnarray}
                  and hence
        \begin{eqnarray}\label{2020.7.12-3} \quad\qquad\,  \;\int_{\Omega}\! \left\{\! \frac{1}{(2\pi)^n}\! \int_{{\mathbb{R}}^n}\!\Big( \frac{1}{2\pi i} \int_{\mathcal{C}} e^{-t\tau}\, \mbox{Tr}\,({\mathbf{q}}_{-2} (x, \xi,\tau) ) d\tau\! \Big) d\xi\!\right\}\! dV\!=\!\Big(\frac{n-1}{(4\pi \mu t)^{n/2}} \! +\!  \frac{1}{(4\pi (2\mu\!+\!\lambda) t)^{n/2}}\!\Big){\mbox{Vol}}(\Omega).\end{eqnarray}

In the above discussion, if we replace $x\in \Omega$ by $\overset{*}{x} \in \Omega^*$, then (\ref{2022.11.9-1}) will become
\begin{eqnarray*}\label{2022.11.9-2}
  \!\! \!\!\!&& \!\!\!\!\!{\mathbf{q}}_{-2} (\overset{*}{x},\xi,\tau) =\frac{1}{\tau- \mu \sum\limits_{l,m=1}^n (g^{lm}(\overset{*}{x}))\xi_l \xi_m }\,{\mathbf{I}}_n +\frac{\mu\!+\!\lambda}{  \big(\tau\!-\! \mu \sum\limits_{l,m=1}^n\! (g^{lm}(\overset{*}{x}))\xi_l \xi_m \big)\big(\tau\!- \!(2\mu\!+\!\lambda)  \sum\limits_{l,m=1}^n \!(g^{lm}(\overset{*}{x}))\xi_l \xi_m \big)}\\
 && \;\;\;\;\times \begin{bmatrix}\! \sum\limits_{r=1}^n \!(g^{1r}(\overset{*}{x})) \xi_r \xi_1 \!&\! \cdots \!&\! \sum\limits_{r=1}^n( g^{1r}(\overset{*}{x}))\xi_r \xi_{n-1}\!& \!\sum\limits_{r=1}^n \!(-g^{1r}(\overset{*}{x})) \xi_r \xi_n\\
\vdots\!& {}\! &\vdots \\
\sum\limits_{r=1}^n (g^{n\!-\!1,r}(\overset{*}{x})) \xi_r \xi_1 \!\!&\!\cdots \!&\! \sum\limits_{r=1}^n (g^{n-1,r}(\overset{*}{x})) \xi_r \xi_{n\!-\!1}\!& \! \sum\limits_{r=1}^n (-g^{n\!-\!1,r}(\overset{*}{x})) \xi_r \xi_n
\\
 \sum\limits_{r=1}^n\! (-g^{nr}(\overset{*}{x})) \xi_r \xi_1 \!&\!\cdots \!& \!\sum\limits_{r=1}^n \!(-g^{nr}(\overset{*}{x})) \xi_r \xi_{n-1}\!& \! \sum\limits_{r=1}^n (g^{nr} (\overset{*}{x}))\xi_r \xi_n\end{bmatrix}\nonumber\end{eqnarray*}
 and   \begin{eqnarray*}\label{2022.11.9-4}  \mbox{Tr}\, \big({\mathbf{q}}_{-2} (\overset{*}{x},\xi,\tau)\big)=
   \frac{n}{\big(\tau \!-\!\mu \!\sum_{l,m=1}^n \!(g^{lm}(\overset{*}{x})) \xi_l\xi_m\big)} \!+\! \frac{(\mu+\lambda)\sum_{l,m=1}^n (g^{lm}(\overset{*}{x})) \xi_l\xi_m}{
  \big(\!\tau\! -\!\mu\! \sum_{l,m=1}^n \!(g^{lm}(\overset{*}{x})) \xi_l\xi_m\!\big)
 \big(\!\tau \!-\!(2\mu\!+\!\lambda)  \!\sum_{l,m=1}^n \!(g^{lm}(\overset{*}{x})) \xi_l\xi_m\!\big)}
.\end{eqnarray*}
This implies that all expressions (\ref{2020.6.6-1})--(\ref{2020.7.12-3}) of the above  trace symbols  have the same form either in $\Omega$ or in $\Omega^*$.

      For given (small) $\epsilon>0$ , denote by $U_\epsilon(\partial \Omega)=\{z\in {\mathcal{M}}\big| \mbox{dist}\, (z, \partial \Omega)<\epsilon\}$ the $\epsilon$-neighborhood of $\partial \Omega$ in $\mathcal{M}$.
             When $x\in \Omega\setminus U_\epsilon (\partial \Omega)$,
we see by taking geodesic normal coordinate system at $x$ that (\ref{2020.7.6-3}) still holds at this $x$. According to (\ref{3.10}) we have
  that
     \begin{eqnarray} \label{3.11}  \mbox{Tr}\,({\mathbf{q}}_{-2}(t, x, \overset{*}{x}))\!\!\! &=\!\!\!&
     \frac{1}{(2\pi)^n} \int_{{\Bbb R}^n} e^{i(x-\overset{*}{x})\cdot\xi}
    \bigg((n-1) e^{-t\mu|\xi|^2} +  e^{-t(2\mu+\lambda)|\xi|^2}\bigg)
         d\xi \nonumber\\
               \!\!\! &=\!\!\!&\frac{n-1}{(4\pi \mu t)^{n/2}}e^{-\frac{|x-\overset{*}{x}|^2}{4t\mu}} +  \frac{1}{(4\pi (2\mu+\lambda) t)^{n/2}}e^{-\frac{|x-\overset{*}{x}|^2}{4t(2\mu+\lambda)}} \quad \, \mbox{for any}\;\, x\in \Omega \setminus U_\epsilon (\partial \Omega),
                  \nonumber\end{eqnarray}
     which exponentially tends to zero as $t\to 0^+$ because $|x-\overset{*}{x}|\ge \epsilon$.  Hence
     \begin{eqnarray} \label{2020.7.6-4}  \int_{\Omega \setminus U_\epsilon (\partial \Omega)} \left(\mbox{Tr}\,({\mathbf{q}}_{-2}(t, x, \overset{*}{x}))\right)\,dV\!\!\! &=\!\!\!&
     O(t^{1-\frac{n}{2}}) \quad \;\mbox{as}\;\; t\to 0^+.
                 \end{eqnarray}

                Secondly, for $l\ge 1$, it can be verified that
                  $\mbox{Tr}\, ({\mathbf{q}}_{-2-l} (x, \xi, \tau))$ is a sum of finitely many terms, each of which  has the following form:
                  $$\frac{r_k(x, \xi)}{(\tau-\mu \sum_{l,m=1}^n g^{lm} \xi_l \xi_m )^s (\tau -(2\mu +\lambda)\sum_{l,m=1}^n g^{lm} \xi_l\xi_m )^j },$$
where $k-2s-2j=-2-l$, and $r_k(x,\xi)$ is the symbol independent of $\tau$ and homogeneous of degree $k$.
 Again we take the geodesic normal coordinate systems center at $x$ (i.e., $g_{jk}(x)=\delta_{jk}$ and $\Gamma_{jk}^l(x)=0$), by applying residue theorem we see  that, for $l\ge 1$ , \begin{eqnarray*}\label{2020.7.12-1} && \frac{1}{(2\pi)^n} \int_{{\mathbb{R}}^n}\Big( \frac{1}{2\pi i} \int_{\mathcal{C}} e^{-t\tau}\, \mbox{Tr}\,({\mathbf{q}}_{-2-l} (x, \xi,\tau) ) d\tau \Big) d\xi=O(t^{l-\frac{n}{2}}) \;\;\mbox{as}\;\, t\to 0^+ \;\;\, \mbox{uniformly for} \;\, x\in \Omega,
 \end{eqnarray*} and
\begin{eqnarray} \label{2020.7.14-1} &\quad\;\;\;\;
 \frac{1}{(2\pi)^n} \!\int_{{\mathbb{R}}^n}\! e^{i(x-\overset{*}{x})\cdot \xi} \Big(\! \frac{1}{2\pi i} \!\int_{\mathcal{C}}\! e^{-t\tau}\, \mbox{Tr}\,({\mathbf{q}}_{-2-l} (x, \xi,\tau) ) d\tau \!\Big) d\xi\!=\!O(t^{l-\frac{n}{2}}\!) \;\,\mbox{as}\;\, t\to 0^+ \,\; \mbox{uniformly for} \;\, x\in \Omega.
      \end{eqnarray}
Therefore
\begin{eqnarray}\label{2020.7.12-10} && \int_{\Omega} \bigg\{\frac{1}{(2\pi)^n} \int_{{\mathbb{R}}^n}\Big( \frac{1}{2\pi i} \int_{\mathcal{C}} e^{-t\tau}\, \sum_{l\ge 1}\mbox{Tr}\,({\mathbf{q}}_{-2-l} (x, \xi,\tau) ) d\tau \Big) d\xi\bigg\} dV =O(t^{1-\frac{n}{2}}) \;\;\mbox{as}\;\, t\to 0^+,
 \end{eqnarray} and
\begin{eqnarray}\label{2020.7.13-5}\;\;\;\;\;
\int_{\Omega} \bigg\{ \frac{1}{(2\pi)^n} \!\int_{{\mathbb{R}}^n}\! e^{i(x-\overset{*}{x})\cdot \xi} \Big( \frac{1}{2\pi i} \!\int_{\mathcal{C}} e^{-t\tau}\, \sum_{l\ge 1}\mbox{Tr}\,({\mathbf{q}}_{-2-l} (x, \xi,\tau) ) d\tau \!\Big) d\xi\bigg\}dV\!=\!O(t^{1-\frac{n}{2}})\;\;\mbox{as}\;\, t\to 0^+.
      \end{eqnarray}
Combining (\ref{2020.7.5-1}), (\ref{2020.7.12-3}) and (\ref{2020.7.12-10}), we have
       \begin{eqnarray} \label{3.11}&& \int_{\Omega}\mbox{Tr}\,({\mathbf{K}}(t, x, x)) \, dV =    \bigg[  \frac{n-1}{(4\pi \mu t)^{n/2}} +  \frac{1}{(4\pi (2\mu+\lambda) t)^{n/2}}\bigg]{\mbox{Vol}}(\Omega)+O(t^{1-\frac{n}{2}})\;\;\mbox{as}\;\; t\to 0^+. \end{eqnarray}

     Finally, we will consider the case of $\int_{\Omega\cap U_\epsilon(\partial \Omega)} \left\{\frac{1}{(2\pi)^n} \int_{{\mathbb{R}}^n} e^{i(x-\overset{*}{x})\cdot \xi} \Big(  \frac{1}{2\pi i} \int_{\mathcal{C}} e^{-t\tau}\,\mbox{Tr}\,({\mathbf{q}}_{-2} (x, \xi,\tau) ) d\tau \Big) d\xi\right\} dV$.
We pick a self-double patch $W$ of $\mathcal{M}$ (such that $W\subset U_\epsilon (\partial \Omega)$) covering a patch $W\cap \partial \Omega$ of $\partial \Omega$  endowed (see the diagram on p.$\,$54 of \cite{MS}) with local coordinates $x$ such that
\begin{figure}[h]
\centering
%\caption{\label{figure1}}
\begin{tikzpicture}[scale=1,line width=0.8]
\clip (-1,1.3) rectangle (6.5,7.4);

\path(70:7) coordinate(A);
\path(30:7) coordinate(B);
\path(70:4.5) coordinate(C);
\path(30:4.5) coordinate(D);

\draw (A) arc (70:30:7);
\draw (C) arc (70:30:4.5);
\draw (A)--(C);
\draw (B)--(D);

\draw[fill=black] (3.177,3.6) circle (1pt);
\draw[fill=black] (4.3,5) circle (1pt);
\draw[fill=black] (3.67,4.32) circle (1pt);

\draw[->,>=stealth] (4.3,5) .. controls (3.5,4.3) and (2.8,3) .. (2.6,2.5);

\node at (68:7.5) {$\Omega^{*}$};
\node at (2.2,2.65) {$\Omega$};
\node at (2.6,2) {$x_n>0$};
\node at (3.35,3.35) {$x$};
\node at (4.5,4.7) {$x^*$};
\node at (-0.5,6.7) {$\partial \Omega$};
\node at (-2,-2){};

\draw (-0.524,6.3709) arc (80.0001:25:8.5);
\end{tikzpicture}
\end{figure}
   $\epsilon>x_n>0$ in $W\cap \Omega$; $\,x_n=0$ on $W\cap \partial \Omega$;
  $\; x_n (\overset{*}{x})=-x_n(x)$; and the positive $x_n$-direction is perpendicular to $\partial \Omega$. This has the effect that (\ref{2021.2.6-3})--(\ref{2021.2.6-4})  and  \begin{eqnarray}
   \label{2021.2.6-5}  \sqrt{|g|/g_{nn}} \; dx_1\cdots dx_{n-1}\!\!\!&\!=\!&\!\!\! \mbox{the element of (Riemannian) surface area on} \,\, \partial \Omega.\end{eqnarray}
    We choose coordinates $x'=(x_1,\cdots, x_{n-1})$ on an open set in $\partial \Omega$ and then
 coordinates $(x', x_{n})$ on a neighborhood in $\bar\Omega$ such that
$x_{n}=0$ on $\partial \Omega$ and $|\nabla x_{n}|=1$ near $\partial \Omega$ while $x_{n}>0$ on $\Omega$ and such that $x'$ is constant
on each geodesic segment in $\bar\Omega$ normal to $\partial \Omega$.
   Then the metric tensor on $\bar \Omega$ has
the form (see \cite{LU} or p.$\,$532 of \cite{Ta2})
\begin{gather} \label{a-1} \big(g_{jk} (x',x_{n}) \big)_{n\times n} =\begin{pmatrix} ( g_{jk} (x',x_{n}))_{(n-1)\times (n-1)}& 0\\
      0& 1 \end{pmatrix}. \end{gather}
       Furthermore, we can take a geodesic normal coordinate system for $(\partial \Omega, g)$ centered at $x_0=0$, with respect to $e_1, \cdots, e_{n-1}$, where  $e_1, \cdots, e_{n-1}$ are the principal curvature vectors. As Riemann showed, one has (see p.$\,$555 of \cite{Ta2})
            \begin{eqnarray} \label{7/14/1}& g_{jk}(x_0)= \delta_{jk}, \; \; \frac{\partial g_{jk}}{\partial x_l}(x_0)
 =0  \;\;  \mbox{for all} \;\; 1\le j,k,l \le n-1,\\
 & -\frac{1}{2}\frac{\partial g_{jk}}{\partial x_n} (x_0) =\kappa_k\delta_{jk}  \;\;  \mbox{for all} \;\; 1\le j,k \le n-1,\nonumber
 \end{eqnarray}
 where $\kappa_1\cdots, \kappa_{n-1}$ are the principal curvatures of $\partial \Omega$ at point $x_0=0$.
   Due to the special geometric normal coordinate system and (\ref{7/14/1})--(\ref{a-1}), we see that for any $x\in \{z\in \Omega\big|\mbox{dist}(z, \partial \Omega)< \epsilon\}$,
   \begin{eqnarray} \label{18-4-5-1} x-\overset{\ast} {x}=(0,\cdots, 0, x_n-(-x_n))=(0,\cdots, 0,2x_n).\end{eqnarray}
              By (\ref{18-4-2-8}), (\ref{7/14/1}), (\ref{3.10}) and (\ref{18-4-5-1}), we find that
      \begin{eqnarray*} \label{18-4-1-1} &&\int_{W\cap \Omega} \bigg\{\frac{1}{(2\pi)^n} \int_{{\Bbb R}^{n}} e^{i\langle x-\overset{*}{x}, \xi\rangle} \Big( \frac{1}{2\pi i} \int_{\mathcal{C}} e^{-t\tau} \,\mbox{Tr}\,\big({\mathbf{q}}_{-2} (x, \xi, \tau)\big) d\tau \Big)  \, d\xi\bigg\} dV \\
           &&   =\! \int_0^\epsilon dx_n\! \int_{W\cap \partial \Omega}\!
 \frac{dx'}{(2\pi)^n}\!\int_{{\Bbb R}^{n}} \!e^{i\langle 0, \xi'\rangle + i2x_n \xi_n} \bigg[\!\frac{1}{2\pi i}\! \int_{\mathcal{C}}\! e^{-t\tau} \!\Big(  \frac{n}{(\tau \!-\!\mu |\xi|^2)}\!+\!
      \frac{(\mu\!+\!\lambda)|\xi|^2}{(\tau \!-\!\mu |\xi|^2) (\tau\!-\!(2\mu\!+\!\lambda)|\xi|^2)}\Big) d\tau\!\bigg] d\xi\\
  && \,=  \int_0^\epsilon dx_n \int_{W\cap \partial \Omega}
 \frac{dx'}{(2\pi)^n}\int_{{\Bbb R}^{n}} e^{i2x_n \xi_n} \bigg((n-1) e^{-t\mu|\xi|^2} +  e^{-t(2\mu+\lambda)|\xi|^2} d\tau\bigg) d\xi\\
 &&\,= \int_0^\epsilon dx_n \int_{W\cap \partial \Omega}  \frac{dx'}{(2\pi)^n}\int_{-\infty}^\infty e^{2ix_n \xi_n}\bigg[ \int_{{\Bbb R}^{n-1}} \bigg( (n-1) e^{-t\mu(|\xi'|^2+\xi_n^2)} +  e^{-t(2\mu+\lambda)(|\xi'|^2+\xi_n^2)} \bigg) d\xi' \bigg]d\xi_n  \nonumber\\
    &&\,= \int_0^\epsilon dx_n \int_{W\cap \partial \Omega}  \frac{1}{(2\pi)^n}\bigg[\int_{-\infty}^\infty e^{2ix_n \xi_n}e^{-t\mu \xi_n^2} \bigg( \int_{{\Bbb R}^{n-1}}  (n-1) e^{-t\mu \sum_{j=1}^{n-1}\xi_j^2}  d\xi' \bigg)d\xi_n\bigg] dx'\nonumber\\
  &&\quad\, + \int_0^\epsilon dx_n \int_{W\cap \partial \Omega}  \frac{1}{(2\pi)^n}\bigg[\int_{-\infty}^\infty e^{2ix_n \xi_n}e^{-t(2\mu+\lambda) \xi_n^2} \bigg( \int_{{\Bbb R}^{n-1}} e^{-t(2\mu+\lambda) \sum_{j=1}^{n-1}\xi_j^2} d\xi' \bigg)d\xi_n\bigg]dx'\nonumber, \end{eqnarray*}
 where  $\xi=(\xi', \xi_n)\in {\Bbb R}^n$, $\xi'=(\xi_1, \cdots, \xi_{n-1})$.
   A direct calculation shows that
   \begin{eqnarray*}
 &&  \frac{1}{(2\pi)^n}\bigg[\int_{-\infty}^\infty e^{2ix_n \xi_n}e^{-t\mu \xi_n^2} \bigg( \int_{{\Bbb R}^{n-1}}  (n-1) e^{-t\mu \sum_{j=1}^{n-1}\xi_j^2}  d\xi' \bigg)d\xi_n\bigg]=\frac{n-1}{(4\pi \mu  t)^{n/2}}\, e^{-\frac{(2x_n)^2}{ 4\mu t} },\\
  &&    \frac{1}{(2\pi)^n}\bigg[\int_{-\infty}^\infty e^{2ix_n \xi_n}e^{-t(2\mu+\lambda) \xi_n^2} \bigg( \int_{{\Bbb R}^{n-1}} e^{-t(2\mu+\lambda) \sum_{j=1}^{n-1}\xi_j^2} d\xi' \bigg)d\xi_n= \frac{1}{(4\pi (2 \mu+\lambda)  t)^{n/2}}\, e^{-\frac{(2x_n)^2}{ 4(2\mu+\lambda) t} }.\end{eqnarray*}
   Hence
    \begin{eqnarray} \label{18-4-1-10.} \\
 && \int_{W\cap \Omega} \bigg\{\frac{1}{(2\pi)^n} \int_{{\Bbb R}^{n}} e^{i\langle x-\overset{*}{x}, \xi\rangle} \Big( \frac{1}{2\pi i} \int_{\mathcal{C}} e^{-t\tau} \,\mbox{Tr}\,\big({\mathbf{q}}_{-2} (x, \xi, \tau)\big) d\tau \Big)  \, d\xi\bigg\} dV\nonumber \\
  &&\;\;\quad = \int_0^\epsilon dx_n \int_{W\cap \partial \Omega} \left[\frac{n-1}{(4\pi \mu  t)^{n/2}}\, e^{-\frac{(2x_n)^2}{ 4\mu t} }+\frac{1}{(4\pi (2 \mu+\lambda)  t)^{n/2}}\, e^{-\frac{(2x_n)^2}{ 4(2\mu+\lambda) t} }\right]dx'
    \nonumber\\
    &&\;\; \quad= \int_0^\infty dx_n \int_{W\cap \partial \Omega} \left[ \frac{n-1}{(4\pi \mu  t)^{n/2}}\, e^{-\frac{(2x_n)^2}{ 4\mu t} }+\frac{1}{(4\pi (2\mu+\lambda)  t)^{n/2}}\, e^{-\frac{(2x_n)^2}{ 4(2\mu+\lambda) t} }\right]dx'\nonumber\\
    && \;\; \quad\;\quad - \int_\epsilon^\infty dx_n \int_{W\cap \partial \Omega} \left[\frac{n-1}{(4\pi \mu  t)^{n/2}}\, e^{-\frac{(2x_n)^2}{ 4\mu t} }+ \frac{1}{(4\pi (2\mu +\lambda) t)^{n/2}}\, e^{-\frac{(2x_n)^2}{ 4(2\mu+\lambda) t} }\right]dx'\nonumber\\
   &&\;\;\quad= \frac{n-1}{4} \cdot \frac{\mbox{Vol}(W\cap \partial \Omega)}{(4\pi \mu t)^{(n-1)/2}} + \frac{1}{4} \cdot \frac{\mbox{Vol}(W\cap \partial \Omega)}{(4\pi (2\mu +\lambda) t)^{(n-1)/2}}\nonumber\\
   &&\;\;\quad \quad \, -   \int_{W\cap \partial \Omega} \left\{\int_\epsilon^\infty\bigg[\frac{n-1}{(4\pi \mu  t)^{n/2}}\, e^{-\frac{(2x_n)^2}{ 4\mu t} }+ \frac{1}{(4\pi (2\mu +\lambda) t)^{n/2}}\, e^{-\frac{(2x_n)^2}{ 4(2\mu+\lambda) t}} \bigg]dx_n\right\}dx'. \nonumber\end{eqnarray}
      It is easy to verify that for any  fixed $\epsilon>0$,
   \begin{eqnarray} \label{18-4-1-11.}\begin{aligned} && \int_\epsilon^\infty \frac{1}{(4\pi \lambda t)^{\frac{n}{2}}} e^{-\frac{(2x_n)^2}{4\mu t}}dx_n = O(t^{1-n/2})\quad \; \, \mbox{as} \;\, t\to 0^+,\qquad \;\; \quad \\
  && \int_\epsilon^\infty \frac{1}{(4\pi (2\mu+\lambda) t)^{\frac{n}{2}}} e^{-\frac{(2x_n)^2}{4(2\mu+\lambda) t}}dx_n = O(t^{1-n/2})
   \quad \; \, \mbox{as} \;\, t\to 0^+.\end{aligned}\end{eqnarray}
   From (\ref{18-4-1-10.}) and (\ref{18-4-1-11.}), we get that
     \begin{eqnarray}\label{18-4-1-6}      && \int_{W\cap \Omega} \bigg\{\frac{1}{(2\pi)^n} \int_{{\Bbb R}^{n}} e^{i\langle x-\overset{*}{x}, \xi\rangle} \Big( \frac{1}{2\pi i} \int_{\mathcal{C}} e^{-t\tau} \,\mbox{Tr}\,\big({\mathbf{q}}_{-2} (x, \xi, \tau)\big) d\tau \Big)  \, d\xi\bigg\} dV=
     \frac{n-1}{4} \cdot \frac{\mbox{Vol}(W\cap \partial \Omega)}{(4\pi \mu t)^{(n-1)/2}}\\
  && \; \quad\,\quad \;\,\quad  + \frac{1}{4} \cdot \frac{\mbox{Vol}(W\cap \partial \Omega)}{(4\pi (2\mu +\lambda) t)^{(n-1)/2}}
 + O(t^{1-n/2}) \quad \; \mbox{as} \, \; t\to 0^+.\nonumber\end{eqnarray}
  For any $x\in \Omega\cap U_\epsilon (\partial \Omega)$, we have
  \begin{eqnarray} \label{18-4-3-1} \mbox{Tr}\,(K(t, x, \overset{*}{x}))
 \!\!\!&\!=&\!\!\!\! \frac{1}{(2\pi)^n}\int_{{\Bbb R}^n}   e^{i\langle x-\overset{*}{x}, \xi\rangle}\big( \frac{1}{2\pi i} \int_{\mathcal{C}} e^{-t\tau} \mbox{Tr}\,\big({\mathbf{q}}_{-2} (x, \xi, \tau)\big) d\tau\big) d\xi\\
  \!\!\!&\!&\!\!\!\!+ \frac{1}{(2\pi)^n}\int_{{\Bbb R}^n}   e^{i\langle x-\overset{*}{x}, \xi\rangle}\big(\sum_{l\ge 1} \frac{1}{2\pi i} \int_{\mathcal{C}} e^{-t\tau} \mbox{Tr}\,\big({\mathbf{q}}_{-2-l} (x, \xi, \tau)\big) d\tau\big) d\xi\nonumber\\
  \!\!\!&\!=&\!\!\!\! \!\frac{1}{(2\pi)^n}\!\int_{{\Bbb R}^n}  \! e^{i\langle x-\overset{*}{x}, \xi\rangle}\big( \frac{1}{2\pi i} \!\int_{\mathcal{C}} e^{-t\tau} \mbox{Tr}\,\big({\mathbf{q}}_{-2} (x, \xi, \tau)\big) d\tau\big) d\xi\!+\!O(t^{1+\frac{n}{2}}) \;\; \mbox{as}\;\, t\to 0^+,\nonumber\end{eqnarray}  where the second equality used (\ref{2020.7.14-1}).
  Combining (\ref{18-4-1-6}) and (\ref{18-4-3-1}), we have
     \begin{eqnarray} \label{3..15}
     &&\int_{W\cap \Omega}\mbox{Tr}\big( \mathbf{K}(t, x, \overset{*}{x})\big)dx=  \frac{n-1}{4} \cdot \frac{\mbox{Vol}(W\cap \partial \Omega)}{(4\pi \mu t)^{(n-1)/2}} \\
  && \quad\,\quad \;\,\quad \;+ \frac{1}{4} \cdot \frac{\mbox{Vol}(W\cap \partial \Omega)}{(4\pi (2\mu +\lambda) t)^{(n-1)/2}}
 + O(t^{1-n/2}) \quad \; \mbox{as} \, \; t\to 0^+.\nonumber\end{eqnarray}
     It follows from  (\ref{c4-23}), (\ref{2020.7.5-2}), (\ref{2020.7.6-4}), (\ref{2020.7.13-5}), (\ref{3.11}) and (\ref{3..15}) that
 \begin{eqnarray} \label{a-4-1-3}  \int_{W\cap \Omega}\mbox{Tr}\big( {\mathbf{K}}^{\mp}(t, x, x) \big) dx\! \!\!\!\! &&\!\!\!\!= \int_{W\cap \Omega}\mbox{Tr}\big({\mathbf{K}}(t, x, x) \big) dx  \mp \int_{W\cap \Omega}\mbox{Tr}\big( {\mathbf{K}}(t, x, \overset{*}{x})\big)dx \\
\! \!\!\!\! &&\!\!\!\!=\bigg[  \frac{n-1}{(4\pi \mu t)^{n/2}}  +  \frac{1}{(4\pi (2\mu+\lambda) t)^{n/2}}\bigg]\mbox{Vol}(W\cap\Omega) \nonumber\\
\! \!\!\!\! &&\!\!\!\!\; \;\;\;\mp \frac{1}{4}\bigg[(n-1) \frac{\mbox{Vol}(W\cap \partial \Omega)}{(4\pi \mu t)^{(n-1)/2}}+ \frac{\mbox{Vol}(W\cap \partial \Omega)}{(4\pi (2\mu +\lambda) t)^{(n-1)/2}}\bigg]\nonumber
  \\   \! \!\!\!\! &&\!\!\!\!\;\;\;\;+O(t^{1-n/2})\quad \; \mbox{as} \, \; t\to 0^+,\nonumber\end{eqnarray}
and hence  (\ref{1-7}) holds.
   $\square$

\vskip 0.28 true cm

\noindent{\bf Remark 4.1.} \ \  i) \  It is very clear that our method and results are still valid in the case of Euclidean space. In other words, if $\Omega$ is a bounded domain in ${\mathbb{R}}^n$, and if $\{\tau_k^-\}$ and $\{\tau_k^{+}\}$  respectively be all the  Navier-Lam\'{e} eigenvalues corresponding to the Navier-Lam\'{e} operator $P\mathbf{u}= -\mu \Delta \mathbf{u} - (\mu+\lambda) \nabla (\nabla\cdot \mathbf{u})$, with the Dirichlet and Neumann boundary condition, then the asymptotic formula (\ref{1-7}) still holds, where the ${\mbox{Vol}}(\Omega)$ and ${\mbox{Vol}}(\partial \Omega)$ are to be replaced by the  $n$-dimensional Euclidean volume $|\Omega|$  and $(n-1)$-dimensional Euclidean volume $|\partial \Omega|$, respectively.

ii) \  \ Note that (see \cite{Avr}, \cite{Avr10}, \cite{Avr11},
\cite{Gil2}, \cite{Mina} or  \cite{Gr})
 $$\int_{\Omega} \mbox{Tr}\big( \mathbf{K}^\mp(t, x,x)) dx = (4\pi t)^{-n/2} \big[a_0+a_1^\mp t^{1/2} + a_2^\mp  t +\cdots +a_m^\mp t^{m/2} +O(t^{(m+1)/2})\big]\;\,\mbox{as}\;\, t\to 0^+.$$  Except for the above obtained $a_0$ and $a_1^\mp$, we can also get all coefficients $a_l^\mp$, $2\le l\le m$, for the asymptotic expansion of the integral of trace of integral kernel for the Navier-Lam\'{e} operator by our new method.

\vskip 0.35 true cm

 Now, we  use the geometric invariants of the Navier-Lam\'{e} spectrum which have been obtained from Theorem 1.1 to finish the proof of Theorem  1.2.

\vskip 0.32 true cm

 \noindent  {\it Proof of Theorem 1.2.} \  By Theorem 1.1, we know that the first two coefficients $a_0$ and $a_1$  of the  asymptotic expansion in (\ref{1-7}) are Navier-Lam\'{e} spectral invariants. From the expressions of $a_0$ and $a_1$, we can further know that $|\Omega|= |B_r|$ and $|\partial \Omega|=|\partial B_r|$. That is, $\frac{|\partial \Omega|}{|\Omega|^{(n-1)/n}} =
 \frac{|\partial B_r|}{|B_r|^{(n-1)/n}}$. Note that for any $r>0$, $\frac{|\partial B_r|}{|B_r|^{(n-1)/n}} =
 \frac{|\partial B_1|}{|B_1|^{(n-1)/n}}$.  According to the classical isoperimetric inequality (which states that for any bounded domain $\Omega\subset {\Bbb R}^n$ with smooth boundary, the following inequality holds:
 \begin{eqnarray*}\frac{|\partial \Omega|}{|\Omega|^{(n-1)/n}} \ge
 \frac{|\partial B_1|}{|B_1|^{(n-1)/n}}.\end{eqnarray*}
 Moreover, equality obtains if and only if $\Omega$ is a ball, see \cite{Ch1} or p.$\,$183 of \cite{CLN}), we immediately get $\Omega=B_r$.  $\;\; \square$

\vskip 0.29 true cm

\noindent{\bf Remark 4.2.} \ \    By applying the Tauberian theorem (see, for example, Theorem 15.3 of p.$\,$30 of \cite{Kor} or p.$\,$446 of \cite{Fel}) for the first term on the right side of (\ref{1-7}) (i.e., $\sum_{k=1}^\infty e^{-t\tau_k^{\mp}}=\int_0^\infty e^{-t\eta} dN^{\mp}(\eta)= \big[ \frac{n-1}{(4\pi \mu t)^{n/2}} +  \frac{1}{(4\pi (2\mu+\lambda) t)^{n/2}}\big]\mbox{Vol}(\Omega)+o(t^{-n/2})$ as $t\to 0^+$), we can easily obtain the Weyl-type law  for the Navier-Lam\'{e} eigenvalues:
\begin{eqnarray} \label{3..2.1} \quad\quad  N^{\mp}(\eta)\!\!\!&=\!&\!\!\max \{k\big| \tau_k^{\mp} \le \eta\} \\
\!&=\!\!&\!\! \frac{\mbox{Vol}(\Omega)}{\Gamma (\frac{n}{2} +1)} \bigg[ \frac{n-1}{(4\pi \mu)^{n/2} }  +\frac{1}{(4\pi (2\mu +\lambda))^{n/2}}\bigg]\eta^{\frac{n}{2}} +o(\eta^{\frac{n}{2}}), \,\quad \; \mbox{as}\;\; \eta \to +\infty.\nonumber\end{eqnarray}.

\vskip 0.22 true cm

\noindent{\bf Remark 4.3.} \ \  Note that  as $\lambda\to -\mu$, the Navier-Lam\'{e} operator reduces to the classical Laplacian.  Therefore
 our results recover all corresponding results for the Laplacian by letting $\lambda+\mu=0$ (cf. \cite{PolS} and \cite{Liu2}).

\vskip 0.22 true cm

\noindent{\bf Remark 4.4.} \ \     Besides the elastic Navier-Lam\'{e} system, by applying our new technique we can also explicitly calculate and obtain the full symbol of the resolvent operator for the thermoelastic system \cite{LiuTan} (it is also an elliptic system). Furthermore, the parabolic trace for this elliptic system can also be explicitly obtained. The thermoelastic system is an important mathematical model, which gives the laws obtained by the deformation of a thermoelastic body and its interior temperature distribution (see Chapter of \cite{LiQi}).

    Generally, for any matrix-valued elliptic operator $L$, if one can explicitly obtain the inverse matrix for the principal symbol of the operator $(\tau I- L)^{-1}$, then all coefficients in the asymptotic expansion of the parabolic trace will  immediately be obtained by our new (algorithm) method.

 \vskip 0.60 true cm

\vskip 0.98 true cm

\centerline {\bf  Acknowledgments}

\vskip 0.39 true cm
  This research was supported by NNSF
of China (11671033/A010802) and NNSF of China (11171023/A010801).

  \vskip 1.68 true cm

\vskip 0.32 true cm

\end{document}